\documentclass[11pt]{article}
\usepackage[tbtags]{amsmath}
\usepackage{amssymb}
\usepackage{amsthm}
\usepackage[misc]{ifsym}
\usepackage{graphicx}
\usepackage{cases}
\usepackage{float}
\numberwithin{equation}{section}   
\normalsize

\title{\bf An Optimal Investment Problem under Correlated Noises: Risk-Sensitive Stochastic Control Approach\thanks{This work is supported by National Key R\&D Program of China (Grant No. 2018YFB1305400) and National Natural Science Foundations of China (Grant No. 11571205, 11831010).}}

\author{\normalsize Le Yang,\thanks{\it Zhongtai Securities Institute for Financial Studies, Shandong University, Jinan 250100, P. R. China E-mail: yangle8886@126.com}\quad
Yueyang Zheng\thanks{\it School of Mathematics, Shandong University, Jinan 250100, P. R. China E-mail: zhengyueyang0106@163.com}\quad
Jingtao Shi\thanks{\it Corresponding author, School of Mathematics, Shandong University, Jinan 250100, P. R. China E-mail: shijingtao@sdu.edu.cn}}

\date{}

\begin{document}

\maketitle

\noindent{\bf Abstract:}\quad This paper is concerned with an optimal investment problem under correlated noises in the financial market, and the expected utility functional is hyperbolic absolute risk aversion (HARA) with the
exponent $\gamma\neq0$. The problem can be reformulated as a risk-sensitive stochastic control problem. A new stochastic maximum principle is obtained first, where the adjoint equations and maximum condition heavily depend on the risk-sensitive parameter and the correlation coefficient. The optimal investment strategy is obtained explicitly in a state feedback form via the solution to a certain Riccati equation, under the risk-seeking case. Numerical simulation and figures are given to illustrate the sensitivity for the optimal investment strategy, with respect to the risk-sensitive parameter and the correlation coefficient.

\vspace{2mm}


\vspace{2mm}


\section{Introduction}

In this paper, we consider an optimal investment problem under correlated noises in the financial market, in which the goal is to maximize the expected utility functional of the wealth. For simplicity, only one risk-free and one risky asset are considered; constraints and transactions costs are ignored. In the traditional Merton's model, the stock price satisfies a geometric Brownian motion. However, we consider a modified model used in Fleming and Sheu \cite{FlemingSheu1999} where the logarithm of the the stock price is subject to an Ornstein-Uhlenbeck type random fluctuation around a deterministic trend, with correlated Brownian noises. We consider a HARA utility functional of the wealth, with exponent $\gamma\neq0$. In Section 2, we reformulated the problem as a risk-sensitive stochastic control problem. The control is the proportion of the wealth invested in the risky asset. We assume that the trend of the logarithm of the stock price is deterministic and linear in the time and the volatility rates are constants. The state is the logarithm of the stock price plus a suitable constant; it satisfies some linear stochastic differential equation (SDE). The problem is then to find a control which maximizes the expectation of an exponential-of-integral cost functional. For more information about the financial model in our paper, please refer to Platen and Rebolledo \cite{PltenRebolledo1996}, Bielecki and Pliska \cite{BieleckiPliska1999}, Pham \cite{Pham2003}, Shi and Wu \cite{ShiWu2012}.

In Section 3, we will set up necessary and sufficient optimal conditions, of the Pontryagin's maximum principle type, for a general risk-sensitive stochastic control problem with correlated Brownian noises. Jacobson \cite{Jacobson1973} has made an initial research of two behaviors with risk-averse and risk-seeking attitudes in the early stage. The dynamic programming principle has been usually the predominant tool to solve the risk-sensitive stochastic control problem. However, like the discussion in this paper, there are several papers having been devoted to the maximum principle. To the best of our knowledge, the pioneer work about risk-sensitive stochastic maximum principle was published by Whittle \cite{Whittle1990}, by the large-deviation theory. In Charalambous and Hibey \cite{CharalambousHibey1996}, a minimum principle for the partial observation risk-sensitive problem is obtained by the measure-valued decomposition and weak control variations. Lim and Zhou \cite{LimZhou2005} obtained a risk-sensitive stochastic maximum principle for the controlled diffusion process with an exponential of integral performance functional, by the relationship between the maximum principle and the dynamic programming principle. Wang and Wu \cite{WangWu2007} derived a risk-sensitive stochastic maximum principle with CRRA's type utility, and applied to solve a portfolio choice problem in the financial market including currency deposit and stock. Then general maximum principle for partially observed risk-sensitive stochastic control problems was proved and applied to finance by Wang and Wu \cite{WangWu2009}, and Huang et al. \cite{HuangLiWang2010}. Shi and Wu \cite{ShiWu2011} discussed the maximum principle for the risk-sensitive stochastic control problem with jump diffusion process, where the control entered both the diffusion and jump terms. Shi and Wu \cite{ShiWu2012} studied a kind of optimal portfolio choice problem in the financial market by the risk-sensitive stochastic maximum principle of \cite{LimZhou2005}, and illustrated the numerical results and figures of optimal investment policies and the sensitivity to the volatility parameter. Djehiche et al. \cite{DjehicheTembineTempone2015} studied the risk-sensitive control problem for system that are non-Markovian and of mean-field type, where the state, the control and the mean of the distribution of state enter the drift term, diffusion term and terminal cost functional. Ma and Liu \cite{MaLiu2016} proved the maximum principle for partially observed risk-sensitive stochastic control problems of mean-field type. Chala \cite{Chala2017} built a stochastic maximum principle for the risk-sensitive control problem for system of backward stochastic differential equation (BSDE), and gave a new method of the transformation of the adjoint process. Sun et al. \cite{SunBrownPamen2018} derived a general stochastic maximum principle for the risk-sensitive optimal control problem of Markov regime-switching jump-diffusion model. Very recently, Moon et al. \cite{MoonDuncanBasar2019} investigated the maximum principle for a two-player risk-sensitive zero-sum differential game. Moon \cite{Moon2019} obtained necessary and sufficient conditions for risk-sensitive stochastic control and differential games with delay. Moon and Basar \cite{MoonBasar2019} considered a risk-sensitive mean-field game via the stochastic maximum principle. For more recent progress for the risk-sensitive stochastic control and differential games with financial applications, please refer to \cite{BensoussanVanSchuppen1985,James1992,FlemingSheu2000,LimZhou2001,NagaiPeng2002,TembineZhuBasar2014,DateGashi2014,HamadeneMu2015,ExarchosTheodorouTsiotras2016,Hata2017,MaLiu2017,MoonBasar2017,XuRM2018,SunPamen2018} and the literatures therein. 

In Section 3, we establish a new risk-sensitive stochastic maximum principle, and the key difference is that we consider the two one-dimension Brownian motions in the state equation and there is a correlation coefficient between them. A new Hamiltonian function, first- and second-order adjoint equations are introduced, which depend on the risk-sensitive parameter and the correlation coefficient. We further prove that the maximum principle is sufficient under some additional convexity/concavity conditions. The detailed proofs will be fulfilled in the Appendix.

In Section 4, we use the new risk-sensitive stochastic maximum principle to solve the optimal investment problem in the financial market, as introduced in Section 2. We give the optimal investment strategy by virtue of the solution to one Riccati equation. The solvability of the Riccati equation is investigated under the risk-seeking case. In Section 5, by virtue of the numerical simulation, some figures are given to illustrate the optimal investment strategy and its sensitivity to the risk-sensitive parameter and the correlation coefficient respectively. Some concluding remarks are given in Section 6.

In this paper, we use $\mathbb{R}^n$ to denote the Euclidean space of $n$-dimensional vectors, $\mathbb{R}^{n\times d}$ to denote the space of $n\times d$ matrices, and $\mathcal{S}^n$ to denote the space of $n\times n$ symmetric matrices. $\langle\cdot,\cdot\rangle$ and $|\cdot|$ are used to denote the scalar product and norm in the Euclidean space, respectively. A $\top$ appearing in the superscript of a matrix, denotes its transpose. $f_x,f_{xx}$ denote the first- and second-order partial derivatives with respect to $x$ for a differentiable function $f$, respectively.

\section{The Model and Problem Formulation}

We consider a financial market, in which investors can choose two kinds of assets: one is risk-free and the other one is risky. The risk-free asset is called bond, whose price $S_{0}(t)$ at time $t$ satisfies the ordinary differential equation (ODE):
\begin{equation}
dS_{0}(t)=rS_{0}(t)dt,
\end{equation}
where $r$ is the constant interest rate. The risky asset is called stock, whose price at time $t$ is $S_{1}(t)$. Let $L(t)=\log S_{1}(t)$ satisfy the Ornstein-Uhlenbeck type process with correlated nosies:
\begin{equation}
dL(t)=c(\bar{L}(t)-L(t))dt+\sigma d\tilde{W}(t)+\bar{\sigma}d\bar{W}(t).
\end{equation}
Here, the price of stock is influenced by two random noises: $\tilde{W}(\cdot)$ and $\bar{W}(\cdot)$, which are both one-dimensional Brownian motions under some probability measure $\tilde{\mathbb{P}}$ on some given measurable space $(\Omega,\mathcal{F})$, and the correlation coefficient between these two Brown motions is a constant $\rho\in[-1,1]$. That is, $d\langle\tilde{W}(t),\bar{W}(t)\rangle=\rho dt$. The stock price volatility rates $\sigma$ and $\bar{\sigma}$ are the nonzero constant, and $c>0$ is some coefficient. $\bar{L}(t)$ is linear in $t$ and is the deterministic log stock price trend. That means that $\bar{L}(t)=mt+\bar{L}_0$, where $m$ and $\bar{L}_0$ are constant.

Let $X(t)$ be the amount of the investor's wealth and $u(t)$ is the proportion of the wealth invested in stock at time $t$, so $(1-u(t))X(t)$ is the amount invested in the bond. In this paper, we require $u(t)\in U=(-\infty,+\infty)$, which means no investment strategy constraints. The wealth dynamics of the investor with the initial state $x_{0}>0$ satisfy
\begin{equation}
\left\{
\begin{aligned}
dX(t)=&\ X(t)\Big\{(1-u(t))rdt+u(t)\big[dL(t)+\frac{1}{2}(\sigma^2+\bar{\sigma}^2+2\rho\sigma\bar{\sigma})dt\big]\Big\},\\
 X(0)=&\ x_0.
\end{aligned}
\right.
\end{equation}
The investor hopes to maximize the expected utility (HARA case) in the terminal time $T>0$:
\begin{equation}\label{HARA}
\tilde{J}(u(\cdot))=\frac{1}{\gamma}\tilde{\mathbb{E}}\big[X(T)^{\gamma}\big],
\end{equation}
by choosing an optimal investment strategy $u(\cdot)\in U$, where $\gamma\neq0$ is the risk-sensitive parameter. The $\tilde{\mathbb{E}}[\cdot]$ is the expectation defined under the probability measure $\tilde{\mathbb{P}}$ while $\mathbb{E}[\cdot]$ is the expectation defined under another probability measure $\mathbb{P}$ which will be introduced later.

Let us now explain the meaning of risk-sensitive by an intuitive argument. Define $\Psi(x)=\frac{1}{\gamma}x^\gamma,\gamma\neq0$. Obviously, $\Psi(\cdot)$ is differentiable at $\tilde{\mathbb{E}}[X(T)]$. Then Taylor's expansion
yields
\begin{equation*}
\begin{aligned}
\Psi(X(T))&\approx\Psi(\tilde{\mathbb{E}}[X(T)])+\Psi'(\tilde{\mathbb{E}}[X(T)])(X(T)-\tilde{\mathbb{E}}[X(T)])\\
          &\quad+\frac{1}{2}\Psi''(\tilde{\mathbb{E}}[X(T)])(X(T)-\tilde{\mathbb{E}}[X(T)])^2.
\end{aligned}
\end{equation*}
If $\Psi(\cdot)$ is strictly concave near $\tilde{\mathbb{E}}[X(T)]$, then $\Psi''(\tilde{\mathbb{E}}[X(T)])=(\gamma-1)\tilde{\mathbb{E}}[X(T)]^{\gamma-2}<0$ which is equivalent to $-\infty<\gamma<1$. This will reduce the overall cost with a large $|X(T)-\tilde{\mathbb{E}}[X(T)]|$ which implies the investor is risk-seeking. Conversely, if $\Psi(\cdot)$ is strictly convex near $\tilde{\mathbb{E}}[X(T)]$, then $\Psi''(\tilde{\mathbb{E}}[X(T)])=(\gamma-1)\tilde{\mathbb{E}}[X(T)]^{\gamma-2}>0$ which is equivalent to $\gamma>1$. This will bring a penalty to the variance term $(X(T)-\tilde{\mathbb{E}}[X(T)])^2$ in the overall cost. In this case, the investor tries to avoid a large deviation of $X(T)$ from its mean $\tilde{\mathbb{E}}[X(T)]$, which implies the investor is risk-averse. Finally, if $\Psi''(\tilde{\mathbb{E}}[X(T)])=(\gamma-1)\tilde{\mathbb{E}}[X(T)]^{\gamma-2}$ is close or equal to 0 which is equivalent to $\gamma\approx1$, then $\tilde{\mathbb{E}}[\Psi(X(T))]\approx\Psi(\tilde{\mathbb{E}}[X(T)])$, in which case the risk-sensitive model
reduces to the risk-neutral one.

To solve the above problem, first we can rewrite the expectation $\tilde{\mathbb{E}}[X(T)^{\gamma}]$ in terms of an expected exponential-of-integral criterion. Applying It\^{o}'s formula to $\ln X(t)^{\gamma}=\gamma\ln X(t)$, we can get
\begin{equation*}
\begin{split}
\tilde{\mathbb{E}}[X(T)^\gamma]=&\ x_0^\gamma\tilde{\mathbb{E}}\Big[\exp\Big\{\Big[\gamma\int_0^T(1-u(t))r+u(t)c\big[\bar{L}(t)-L(t)+\frac{1}{2}(u(t)-u(t)^2)\big]\\
                                &\qquad\times(\sigma^2+\bar{\sigma}^2+2\rho\sigma\bar{\sigma})\Big]dt+\gamma\int_0^Tu(t)\sigma d\tilde{W}(t)+\gamma\int_0^Tu(t)\bar{\sigma} d\bar{W}(t)\Big\}\Big].
\end{split}
\end{equation*}
We eliminate the stochastic integral term by using the Girsanov transformation:
\begin{equation*}
\begin{split}
\frac{d\mathbb{P}}{d\tilde{\mathbb{P}}}=\exp\Big\{\gamma\int_0^Tu(t)\sigma d\tilde{W}(t)+\gamma\int_0^Tu(t)\bar{\sigma}d\bar{W}(t)-\frac{1}{2}\gamma^2\int_0^T(\sigma^2+\bar{\sigma}^2)u^2(t)dt\Big\}.
\end{split}
\end{equation*}
In order to make the change of probability measure argument valid (see Liptser and Shiryayev \cite{LiptserShiryayev1977}), we assume that existing positive constants $\beta$ and $C$ satisfy the following inequality
\begin{equation}\label{condition1}
\tilde{\mathbb{E}}\Big[\exp\Big\{\beta[(\gamma\sigma u(t))^{2}+(\gamma\bar{\sigma} u(t))^2]\Big\}\Big]\leq C.
\end{equation}
Then
\begin{equation}\label{Lt}
\begin{split}
dL(t)=\big[c(\bar{L}(t)-L(t))+\gamma\sigma^2u(t)+\gamma\bar{\sigma}^2u(t)\big]dt+\sigma dW_1(t)+\bar{\sigma}dW_2(t),
\end{split}
\end{equation}
where $W_1(\cdot)$ and $W_2(\cdot)$ are Brownian motions under the probability measure $\mathbb{P}$ with correlation coefficient still being $\rho$, and
\begin{equation*}
\tilde{\mathbb{E}}[X(T)^{\gamma}]=x_0^\gamma\mathbb{E}\Big[\exp\Big\{\gamma\int_0^Th(x(t),u(t))dt\Big\}\Big],
\end{equation*}
where
\begin{equation}\label{cost function}
\begin{split}
h(x(t),u(t))&=\frac{1}{2}u(t)^2\big[(\gamma-1)(\sigma^2+\bar{\sigma}^2)-2\rho\sigma\bar{\sigma}\big]+r[1-u(t)]\\
            &\quad+u(t)\big[c(\bar{L}(t)-L(t))+\frac{1}{2}(\sigma^2+\bar{\sigma}^2+2\rho\sigma\bar{\sigma})\big].
\end{split}
\end{equation}
Here $L(\cdot)$ is the state and $u(\cdot)$ is the control. We required $u(\cdot)$ is $\mathcal{F}_t$-adapted to the filtration $\mathcal{F}_t=\sigma\{W_1(s),W_2(s);0\leq s\leq t\}$ and (\ref{condition1}) holds. We call such $u(\cdot)$ an admissible control and denote by $\mathcal{U}_{ad}$ the admissible control set.

In order to simplify the problem, we replace $L(t)$ by an equivalent state variable
\begin{equation}\label{state}
x(t):=L(t)-\bar{L}(t)+c^{-1}m.
\end{equation}
Then, by (\ref{Lt}), we have
\begin{equation}\label{state equation}
\left\{
\begin{aligned}
dx(t)=&\ \big[-cx(t)+\gamma(\sigma^2+\bar{\sigma}^2)u(t)\big]dt+\sigma dW_1(t)+\bar{\sigma}dW_2(t),\\
 x(0)=&\ c^{-1}m.
\end{aligned}
\right.
\end{equation}
So the expected HARA utility (\ref{HARA}) is reduced to the expected exponential-of-integral form
\begin{equation}\label{cost functional}
J(u(\cdot))=\mathbb{E}\Big[\exp\Big\{\gamma\int_0^Th(x(t),u(t))dt\Big\}\Big],
\end{equation}
where $h(x,u)$ is defined by (\ref{cost function}). Hence we need to maximize (\ref{cost functional}) by choosing $u(\cdot)$ over $\mathcal{U}_{ad}$.

In order to solve this problem, in the next section, we consider a general risk-sensitive stochastic control problem with two Brownian noises whose correlation coefficient is $\rho$ first. Then we apply the theoretic results to solve this problem in Section 5.

\section{A General Risk-Sensitive Stochastic Control problem}

In this section, we consider a general risk-sensitive stochastic control problem with correlated Brownian noises and give the necessary and sufficient conditions for the optimality.

Let the time duration $T>0$. For any $(s,x)\in[0,T]\times\mathbb{R}^n$, we consider the following SDE:
\begin{equation}\label{SDE}
\left\{
\begin{aligned}
dx(t)=&\ b(t,x(t),u(t))dt+\sigma_1(t,x(t),u(t))dW_1(t)+\sigma_2(t,x(t),u(t))dW_2(t),\\
 x(s)=&\ x.
\end{aligned}
\right.
\end{equation}
We will work in the weak formulation. For any $s\in[0,T]$, the class of admissible controls $\mathcal{U}[s,T]$ is the set of all six-tuple $(\Omega,\mathcal{F},\mathbb{P},W_1(\cdot),W_2(\cdot),u(\cdot))$ satisfying the following conditions.

1) $(\Omega,\mathcal{F},\mathbb{P})$ is a complete probability space.

2) $W_1(\cdot),W_2(\cdot)$ are both one-dimensional Brownian motions defined on $(\Omega,\mathcal{F},\mathbb{P})$ with the correlation coefficient being $\rho\in[-1,1]$, and $\mathcal{F}^{s}_{t}$ is $\sigma\{W_1(r),W_2(r);s\leq r\leq t\}$ augmented by all the $\mathbb{P}$-${\rm null}$ set in $\mathcal{F}$.

3) $u:[s,T]\times\Omega\rightarrow\mathbb{R}^k$ is an $\{\mathcal{F}^s_t\}_{t\geq s}$ adapted process on $(\Omega,\mathcal{F},\mathbb{P})$.

4) Under $u(\cdot)$, for any $x\in\mathbb{R}^n$, (\ref{SDE}) admits a unique adapted solution $x(\cdot)$ on $(\Omega,\mathcal{F},\mathbb{P};\{\mathcal{F}^s_t\}_{t\geq 0})$.

If there are no ambiguity, we will only write $u(\cdot)\in\mathcal{U}[s,T]$ instead of the entire six-tuple $(\Omega,\mathcal{F},\mathbb{P},W_1(\cdot),W_2(\cdot),u(\cdot))$. If $x(\cdot)$ is the unique solution to (\ref{SDE}) associated with the input $u\in\mathcal{U}[s,T]$, we refer to $(x(\cdot),u(\cdot))$ as an admissible pair. In this paper, $L^2_\mathcal{F}(s,T;\mathbb{R}^n)$ will denote the set of $\mathbb{R}^n$-valued, $\{\mathcal{F}^s_{t}\}_{t\geq s}$-adapted, square integrable processes on $[s,T]$, and $L^\infty_\mathcal{F}(s,T;\mathbb{R}^n)$ will denote the set of $\mathbb{R}^n$-valued, $\{\mathcal{F}^s_t\}_{t\geq s}$-adapted, essentially bounded processes on $[s,T]$.

The cost functional $J^{\theta}(s,x,u(\cdot))$ associated with the initial condition $(s,x)\in[0,T]\times\mathbb{R}^{n}$ and $u(\cdot)\in\mathcal{U}[s,T]$ is given by
\begin{equation}\label{cost}
J^\theta(s,x;u(\cdot))=\mathbb{E}\Big[\exp\Big\{\theta\big[g(x(T))+\int_s^Tf(t,x(t),u(t))dt\big]\Big\}\Big],
\end{equation}
where $\theta>0$, the risk-sensitive parameter, is a given fixed constant. The risk-sensitive stochastic control problem associated with (\ref{SDE}) and (\ref{cost}) is defined as follows:
\begin{align}\label{problem}
\begin{cases}
{\rm Maximize}: &J^\theta(s,x;u(\cdot)),\\
{\rm subject\ to}:&
u(\cdot)\in\mathcal{U}[s,T],\ (x(\cdot),u(\cdot))\ \mbox{satisfies}\ (\ref{cost}).
\end{cases}
\end{align}
The value function $v^\theta:[0,T]\times \mathbb{R}^{n}\rightarrow \mathbb{R}$ associated with (7) is defined by
\begin{equation}\label{value function}
v^\theta(s,x):=\sup_{u\in\mathcal{U}[s,T]}J^\theta(s,x;u(\cdot)).
\end{equation}

Let we introduce the following assumptions:

\vspace{1mm}

{\bf(B1)} The map $b:[0,T]\times\mathbb{R}^n\times\mathbb{R}^k\rightarrow\mathbb{R}^n$, $\sigma_i:[0,T]\times\mathbb{R}^n\times\mathbb{R}^k\rightarrow\mathbb{R}^n$, where $i=1,2$, $f:[0,T]\times\mathbb{R}^n\times\mathbb{R}^k\rightarrow\mathbb{R}$ and $g:\mathbb{R}^n\rightarrow\mathbb{R}$ are measurable, and there exists a constant $L>0$ and a modulus of continuity $\bar{\omega}:[0,\infty)\rightarrow[0,\infty)$ such that for $\varphi(t,x,u)=b(t,x,u),\sigma_i(t,x,u),f(t,x,u),g(x)$,
\begin{align*}
\begin{cases}
|\varphi(t,x,u)-\varphi(t,y,u)|\leq L|x-y|+\bar{\omega}(d(u,v)),\ \forall t\in[0,T],\  x,y\in \mathbb{R}^n,\ u,v\in\mathbb{R}^k,\\
|\varphi(t,0,u)|\leq L, \ u\in\mathbb{R}^k.
\end{cases}
\end{align*}
Also $f$ and $g$ are uniformly bounded.

\vspace{1mm}

{\bf(B2)} $b$, $f$ are twice differentiable in $x$, and there exists a modulus of continuity $\bar{\omega}:[0,\infty]\rightarrow[0,\infty]$ such that for $\varphi(t,x,u)=b(t,x,u),\sigma_i(t,x,u),f(t,x,u),g(x)$,
\begin{align*}
\begin{cases}
|\varphi_x(t,x,u)-\varphi_x(t,y,u)|\leq L|x-y|+\bar{\omega}(d(u,v)),\\
|\varphi_{xx}(t,x,u)-\varphi_{xx}(t,y,u)|\leq\bar{\omega}(|x-y|+d(u,v)),\ \forall t\in[0,T],\ x,y\in \mathbb{R}^n,\ u,v\in\mathbb{R}^k.\\
\end{cases}
\end{align*}

{\bf(B3)} $v^\theta\in C^{1,3}([0,T]\times\mathbb{R}^n).$

\vspace{1mm}

{\bf(B4)} $U$ is a convex subset of $\mathbb{R}^k$. The map $b,\sigma,f$ are locally Lipschitz in $u$, and their derivatives in $x$ are continuous in $(x,u)$.

Let $(x(\cdot),u(\cdot))$ be an admissible pair for (\ref{problem}). We introduce the first- and second-order adjoint variables $(\bar{p}(\cdot),\bar{q}_1(\cdot),\bar{q}_2(\cdot))$, $(\bar{P}(\cdot),\bar{Q}_1(\cdot),\bar{Q}_2(\cdot))$, which are the solution to the following BSDEs, respectively:
\begin{equation}\label{BSDE1}
\left\{
\begin{aligned}
d\bar{p}(t)=&-\Big\{{\bar{b}_x(t)}^\top\bar{p}(t)-\bar{f}_x(t)^\top+\bar{\sigma}_{1x}(t)^\top\bar{q}_1(t)-\theta\bar{p}(t)^\top\bar{\sigma}_1(t)\bar{\sigma}_{1x}(t)^\top\bar{p}(t)+\bar{\sigma}_{2x}(t)^\top\bar{q}_1(t)\\
            &\quad-\theta\bar{p}(t)^\top\bar{\sigma}_2(t)\bar{\sigma}_{2x}(t)^\top\bar{p}(t)-\theta\big[\bar{p}(t)^\top\bar{\sigma}_1(t)\bar{q}_1(t)+\rho\bar{p}(t)^\top\bar{\sigma}_2(t)\bar{q}_1(t)\\
            &\quad+\rho\bar{p}(t)^\top\bar{\sigma}_1(t)\bar{q}_2(t)+\bar{p}(t)^\top\bar{\sigma}_2(t)\bar{q}_2(t)\big]\Big\}dt+\bar{q}_1(t)dW_1(t)+\bar{q}_2(t)dW_2(t),\\
 \bar{p}(T)=&-g_x(\bar{x}(T)),
\end{aligned}
\right.
\end{equation}
\begin{equation}\label{BSDE2}
\left\{
\begin{aligned}
d\bar{P}(t)=&-\Big\{\bar{b}_x(t)^\top\bar{P}(t)+\bar{P}(t)\bar{b}_x(t)+\bar{\sigma}_{1x}(t)^\top(\bar{P}(t)-\theta\bar{p}(t)\bar{p}(t)^\top)\bar{\sigma}_{1x}(t)\\
            &\quad+\bar{\sigma}_{2x}(t)^\top(\bar{P}(t)-\theta\bar{p}(t)\bar{p}(t)^\top)\bar{\sigma}_{2x}(t)+\bar{\sigma}_{1x}(t)^\top\big[\bar{Q}_1(t)-\theta\bar{p}(t){\bar{q}_1(t)}^\top\\
            &\quad-\theta{\bar{p}(t)}^\top\bar{\sigma}_1(t)\bar{P}(t)\big]+\bar{\sigma}_{2x}(t)^\top\big[\bar{Q}_2(t)-\theta\bar{p}(t)\bar{q}_2(t)^\top-\theta\bar{p}(t)^\top\bar{\sigma}_2(t)\bar{P}(t)\big]\\
            &\quad+\big[\bar{Q}_1(t)-\theta\bar{q}_1(t)\bar{p}(t)^\top-\theta\bar{p}(t)^\top\bar{\sigma}_1(t)\bar{P}(t)\big]\bar{\sigma}_{1x}(t)\\
            &\quad+\big[\bar{Q}_2(t)-\theta\bar{q}_2(t)\bar{p}(t)^\top-\theta\bar{p}(t)^\top\bar{\sigma}_2(t)\bar{P}(t)\big]\bar{\sigma}_{2x}(t)\\
            &\quad+\bar{H}_{xx}^\theta(t,\bar{x}(t),\bar{u}(t),\bar{p}(t),\bar{q}_1(t),\bar{q}_2(t))-\theta\bar{p}(t)^\top\bar{\sigma}_1(t)\bar{Q}_1(t)-\theta\bar{p}(t)^\top\bar{\sigma}_2(t)\bar{Q}_2(t)\\
            &\quad-\theta\rho\bar{p}(t)^\top\bar{\sigma}_1(t)\bar{Q}_2(t)-\theta\rho\bar{p}(t)^\top\bar{\sigma}_2(t)\bar{Q}_1(t)+\theta^2\rho\bar{p}(t)^\top\bar{\sigma}_1(t)\bar{p}(t)\bar{q}_2(t)^\top\\
            &\quad+\theta^2\rho\bar{p}(t)^\top\bar{\sigma}_1(t)\bar{q}_2(t)\bar{p}(t)^\top+\theta^2\rho\bar{p}(t)^\top\bar{\sigma}_2(t)\bar{p}(t)\bar{q}_1(t)^\top
             +\theta^2\rho\bar{p}(t)^\top\bar{\sigma}_2(t)\bar{q}_1(t)\bar{p}(t)^\top\\
            &\quad-\theta^2\rho\bar{p}(t)\bar{q}_1(t)^\top\bar{\sigma}_2(t)^\top\bar{p}(t)-\theta^2\rho\bar{p}(t)\bar{q}_2(t)^\top\bar{\sigma}_1(t)^\top\bar{p}(t)
             -\theta^2\rho\bar{p}(t)\bar{\sigma}_2(t)^\top\bar{q}_1(t)^\top\bar{p}(t)\\
            &\quad-\theta^2\rho\bar{p}(t)\bar{\sigma}_1(t)^\top\bar{q}_2(t)^\top\bar{p}(t)-\theta(\bar{q}_1(t)\bar{q}_1(t)^\top+\rho\bar{q}_1(t)\bar{q}_2(t)^\top\\
            &\quad+\rho\bar{q}_2(t)\bar{q}_1(t)^\top+\bar{q}_2(t)\bar{q}_2(t)^\top)\Big\}dt+\bar{Q}_1(t)dW_1(t)+\bar{Q}_2(t)dW_2(t),\\
 \bar{P}(T)=&-g_{xx}(\bar{x}(T)),
\end{aligned}
\right.
\end{equation}
where $\bar{b}_x(t):=b_x(t,\bar{x}(t),\bar{u}(t))$ $(\bar{f}_x(t),\bar{\sigma}_{ix}(t),i=1,2$ has the similar interpretations). And the Hamiltonian function $\bar{H}^\theta:\mathbb{R}\times\mathbb{R}^n\times\mathbb{R}^k\times \mathbb{R}^n\times\mathbb{R}^n\times\mathbb{R}^n\rightarrow\mathbb{R}$ is given by
\begin{equation}\label{H1}
\begin{aligned}
\bar{H}^\theta(t,x,u,p,q_1,q_2)=&\left\langle p,b(t,x,u)\right\rangle-f(t,x,u)+\sigma_1(t,x,u)^\top[q_1-\theta pp^\top\sigma_1(t,\bar{x},\bar{u})]\\
                                &+\sigma_{2}(t,x,u)^\top[q_2-\theta pp^\top\sigma_2(t,\bar{x},\bar{u})].
\end{aligned}
\end{equation}
Though (\ref{BSDE1}) is a nonlinear equation, it will be shown that our assumption are sufficient to guarantee the existence of unique solutions $(\bar{p}(\cdot),\bar{q}_1(\cdot),\bar{q}_2(\cdot))\in L^2_\mathcal{F}(s,T;\mathbb{R}^n)\times L^2_\mathcal{F}(s,T;\mathbb{R}^n)\times L^2_\mathcal{F}(s,T;\mathbb{R}^n)$ and $(\bar{P}(\cdot),\bar{Q}_1(\cdot),\bar{Q}_2(\cdot))\in L^2_\mathcal{F}(s,T;\mathbb{R}^{n\times n})\times L^{2}_\mathcal{F}(s,T;\mathbb{R}^{n\times n})\times L^2_\mathcal{F}(s,T;\mathbb{R}^{n\times n})$ to (\ref{BSDE1}) and (\ref{BSDE2}), respectively.

\vspace{1mm}

{\bf Theorem 3.1 (Risk-Sensitive Maximum Principle)} {\it Suppose that {\bf(B1)-(B3)} hold. Let $(\bar{x}(\cdot),\bar{u}(\cdot))$ be an optimal pair for the risk-sensitive stochastic control problem (\ref{problem}). Then, there are unique solutions $(\bar{p}(\cdot),\bar{q}_1(\cdot),\bar{q}_2(\cdot))$ and $(\bar{P}(\cdot),\bar{Q}_1(\cdot),\bar{Q}_2(\cdot))$ to the first-order and the second-order adjoint equations (\ref{BSDE1}) and (\ref{BSDE2}) respectively, such that
\begin{equation}\label{mp1}
\begin{split}
&\bar{H}^\theta(t,\bar{x}(t),\bar{u}(t),\bar{p}(t),\bar{q}_1(t),\bar{q}_2(t))-\bar{H}^\theta(t,\bar{x}(t),u(t),\bar{p}(t),\bar{q}_1(t),\bar{q}_2(t))\\
&-\frac{1}{2}[\sigma_1(t,\bar{x}(t),\bar{u}(t))-\sigma_1(t,\bar{x}(t),u(t))]^\top[\bar{P}(t)-\theta\bar{p}(t)\bar{p}(t)^\top][\sigma_1(t,\bar{x}(t),\bar{u}(t))-\sigma_1(t,\bar{x}(t),u(t))]\\
&-\frac{1}{2}[\sigma_2(t,\bar{x}(t),\bar{u}(t))-\sigma_2(t,\bar{x}(t),u(t))]^\top[\bar{P}(t)-\theta\bar{p}(t)\bar{p}(t)^\top][\sigma_2(t,\bar{x}(t),\bar{u}(t))-\sigma_2(t,\bar{x}(t),u(t))]\\
&\leq0,\ a.e.t\in[s,T],\ \mathbb{P}\mbox{-}{\rm a.s}.
\end{split}
\end{equation}
or equivalently,
\begin{equation}\label{mp2}
\bar{\mathcal{H}}^\theta(t,\bar{x}(t),\bar{u}(t))=\min_{u\in\mathbb{R}^k}\bar{\mathcal{H}}^\theta(t,\bar{x}(t),u),\ a.e. t\in[s,T],\ \mathbb{P}\mbox{-}{\rm a.s.},
\end{equation}
where we define $\mathcal{H}$-function $\bar{\mathcal{H}}^\theta:\mathbb{R}\times\mathbb{R}^n\times\mathbb{R}^k\rightarrow\mathbb{R}$ as
\begin{equation}\label{H2}
\begin{split}
\bar{\mathcal{H}}^\theta(t,x,u):=&\ \big\langle\bar{p}(t),b(t,x,u)\big\rangle-f(t,x,u)+\frac{1}{2}\sigma_1(t,x,u)^\top\big[\bar{P}(t)-\theta\bar{p}(t)\bar{p}(t)^\top\big]\sigma_1(t,x,u)\\
                                 &+\frac{1}{2}\sigma_2(t,x,u)^\top\big[\bar{P}(t)-\theta\bar{p}(t)\bar{p}(t)^\top\big]\sigma_2(t,x,u)+\sigma_1(t,x,u)^\top\big[\bar{q}_1(t)\\
                                 &-\bar{P}(t)\sigma_1(t,\bar{x}(t),\bar{u}(t))\big]+\sigma_2(t,x,u)^\top\big[\bar{q}_2(t)-\bar{P}(t)\sigma_2(t,\bar{x}(t),\bar{u}(t))\big].
\end{split}
\end{equation}}

Sufficient conditions for optimality of the pair $(\bar{x}(\cdot),\bar{u}(\cdot))$ are as follows.

\vspace{1mm}

{\bf Theorem 3.2 (Sufficient Conditions for Optimality)} {\it Suppose that {\bf(B1)-(B4)} hold. Let $(\bar{x}(\cdot),\bar{u}(\cdot))$ be an admissible pair, and $(\bar{p}(\cdot),\bar{q}_1(\cdot),\bar{q}_2(\cdot))$, $(\bar{P}(\cdot),\bar{Q}_1(\cdot),\bar{Q}_2(\cdot))$ be the first- and second-order adjoint variables, respectively. Suppose $g(\cdot)$ is convex, $\bar{H}^\theta(t,\cdot,\cdot,\bar{p}(t),\bar{q}_1(t),\bar{q}_2(t))$ is concave for all $t\in[0,T]$ almost surely and (\ref{mp2}) holds. Then $(\bar{x}(\cdot),\bar{u}(\cdot))$ is an optimal pair for problem (\ref{problem}).}

\vspace{1mm}

Proof of the two theorems will be given in the Appendix.

\section{Optimal Investment Problem: Explicit Solution in the Risk-Seeking Case}

In this section, we will apply the results we got in the previous section and solve the problem in Section 2.

Let $(\bar{x}(\cdot),\bar{u}(\cdot))$ be an optimal pair. First of all, the first-order adjoint equation is as follows:
\begin{equation}\label{adjoint equation}
\left\{
\begin{aligned}
d\bar{p}(t)=&\Big\{\gamma\bar{p}(t)\big[(\bar{q}_1(t)+\bar{q}_2(t)\rho)\sigma+(\bar{q}_2(t)+\bar{q}_1(t)\rho)\bar{\sigma}\big]+c\bar{p}(t)-c\bar{u}(t)\Big\}dt\\
            &+\bar{q}_1(t)dW_{1}(t)+\bar{q}_2(t)dW_2(t),\\
\bar{p}(T)=&\ 0.
\end{aligned}
\right.
\end{equation}
Noting that in our state dynamics (\ref{state equation}), the diffusion term is control independent, so the second-order adjoint variables disappear automatically and $\mathcal{H}$-function gets the form
\begin{equation}
\begin{split}
\bar{\mathcal{H}}^\theta(t,x,u)=&\big[-cx(t)+\gamma(\sigma^2+\bar{\sigma}^2)u(t)\big]\bar{p}(t)-\frac{1}{2}(\sigma^2+\bar{\sigma}^2)\gamma\bar{p}^2(t)\\
                                &-\frac{u^{2}(t)}2\big[(\gamma-1)(\sigma^2+\bar{\sigma}^2)-2\rho\sigma\bar{\sigma}\big]\\
                                &-\Big[\frac{\sigma^2+\bar{\sigma}^2+2\rho\sigma\bar{\sigma}}{2}+m-r-cx(t)\Big]u(t)+\sigma\bar{q}_1(t)+\bar{\sigma}\bar{q}_2(t)-r.
\end{split}
\end{equation}

By the minimum condition (\ref{mp2}), we obtain
\begin{equation}\label{optimal control}
\begin{split}
\bar{u}(t)=&\ \frac{\gamma(\sigma^2+\bar{\sigma}^2)}{(\gamma-1)(\sigma^2+\bar{\sigma}^2)-2\rho\sigma\bar{\sigma}}\bar{p}(t)+\frac{c}{(\gamma-1)(\sigma^2+\bar{\sigma}^2)-2\rho\sigma\bar{\sigma}}\bar{x}(t)\\
           &-\frac{\frac{1}{2}(\sigma^2+\bar{\sigma}^2+2\rho\sigma\bar{\sigma})+m-r}{(\gamma-1)(\sigma^2+\bar{\sigma}^2)-2\rho\sigma\bar{\sigma}}.
\end{split}
\end{equation}
Given $x(0)=c^{-1}m$ and substituting (\ref{optimal control}) into (\ref{state equation}) and (\ref{adjoint equation}), gives the following forward-backward SDE (FBSDE):
\begin{equation}\label{FBSDE}
\left\{
\begin{aligned}
\begin{split}
d\bar{x}(t)=&\bigg\{\frac{(\sigma^2+\bar{\sigma}^2+2\rho\sigma\bar{\sigma})c}{(\gamma-1)(\sigma^2+\bar{\sigma}^2)-2\rho\sigma\bar{\sigma}}\bar{x}(t)
             +\frac{\gamma^2(\sigma^2+\bar{\sigma}^2)^2}{(\gamma-1)(\sigma^2+\bar{\sigma}^2)-2\rho\sigma\bar{\sigma}}\bar{p}(t)\\
            &\quad-\frac{\gamma(\sigma^2+\bar{\sigma}^2)[\frac{1}{2}(\sigma^2+\bar{\sigma}^2+2\rho\sigma\bar{\sigma})+m-r]}{(\gamma-1)(\sigma^2+\bar{\sigma}^2)-2\rho\sigma\bar{\sigma}}\bigg\}dt
             +\sigma dW_1(t)+\bar{\sigma}dW_2(t),\\
d\bar{p}(t)=&\bigg\{\frac{-c(\sigma^2+\bar{\sigma}^2+2\rho\sigma\bar{\sigma})}{(\gamma-1)(\sigma^2+\bar{\sigma}^2)-2\rho\sigma\bar{\sigma}}\bar{p}(t)
             +\frac{\frac{c}{2}(\sigma^2+\bar{\sigma}^2+2\rho\sigma\bar{\sigma})+mc-rc}{(\gamma-1)(\sigma^2+\bar{\sigma}^2)-2\rho\sigma\bar{\sigma}}\\
            &\quad+\gamma\bar{p}(t)\big[(\bar{q}_1(t)+\bar{q}_2(t)\rho)\sigma+(\bar{q}_2(t)+\bar{q}_1(t)\rho)\bar{\sigma}\big]\\
            &\quad-\frac{c^2}{(\gamma-1)(\sigma^2+\bar{\sigma}^2)-2\rho\sigma\bar{\sigma}}\bar{x}(t)\bigg\}dt+\bar{q}_1(t)dW_1(t)+\bar{q}_2(t)dW_2(t),\\
 \bar{x}(0)=&\ x_0,\quad\bar{p}(T)=0.
\end{split}
\end{aligned}
\right.
\end{equation}

As in Yong and Zhou \cite{YongZhou1999}, we conjecture the solution to (\ref{FBSDE}) is related by
\begin{equation}\label{relation}
\bar{p}(t)=-Q(t)\bar{x}(t)-\varphi(t),
\end{equation}
where $Q(\cdot)$ and $\varphi(\cdot)$ are some deterministic differentiable functions. Applying It\^{o}'s formula to (\ref{relation}), it gives
\begin{equation}\label{bar p1}
\begin{split}
d\bar{p}(t)=&\bigg\{\bigg[-\dot{Q}(t)-\frac{c(\sigma^2+\bar{\sigma}^2+2\rho\sigma\bar{\sigma})}{(\gamma-1)(\sigma^2+\bar{\sigma}^2)-2\rho\sigma\bar{\sigma}}Q(t)
             +\frac{\gamma^2(\sigma^2+\bar{\sigma}^2)^2}{(\gamma-1)(\sigma^2+\bar{\sigma}^2)-2\rho\sigma\bar{\sigma}}Q^2(t)\bigg]\bar{x}(t)\\
            &-\dot{\varphi}(t)+\frac{\gamma^2(\sigma^2+\bar{\sigma}^2)^2}{(\gamma-1)(\sigma^2+\bar{\sigma}^2)-2\rho\sigma\bar{\sigma}}Q(t)\varphi(t)\\
            &+Q(t)\frac{\gamma(\sigma^2+\bar{\sigma}^2)[\frac{1}{2}(\sigma^2+\bar{\sigma}^2+2\rho\sigma\bar{\sigma})+m-r]}{(\gamma-1)(\sigma^2+\bar{\sigma}^2)-2\rho\sigma\bar{\sigma}}\bigg\}dt
             -\sigma Q(t)dW_1(t)-\bar{\sigma}Q(t)dW_2(t).
\end{split}
\end{equation}
On the other hand, substituting (\ref{relation}) into (\ref{adjoint equation}), we get
\begin{equation}\label{bar p2}
\begin{split}
d\bar{p}(t)=&\bigg\{\bigg[\frac{c(\sigma^2+\bar{\sigma}^2+2\rho\sigma\bar{\sigma})}{(\gamma-1)(\sigma^2+\bar{\sigma}^2)-2\rho\sigma\bar{\sigma}}Q(t)
             -\frac{c^2}{(\gamma-1)(\sigma^2+\bar{\sigma}^2)-2\rho\sigma\bar{\sigma}}-\gamma[(\bar{q}_1(t)\\
             &+\bar{q}_2(t)\rho)\sigma+(\bar{q}_2(t)+\bar{q}_1(t)\rho)\bar{\sigma}]Q(t)\bigg]\bar{x}(t)
              +\frac{c(\sigma^2+\bar{\sigma}^2+2\rho\sigma\bar{\sigma})}{(\gamma-1)(\sigma^2+\bar{\sigma}^2)-2\rho\sigma\bar{\sigma}}\varphi(t)\\
             &+\frac{c[\frac{1}{2}(\sigma^2+\bar{\sigma}^2+2\rho\sigma\bar{\sigma})+m-r]}{(\gamma-1)(\sigma^2+\bar{\sigma}^2)-2\rho\sigma\bar{\sigma}}-\gamma\big[(\bar{q}_1(t)+\bar{q}_2(t)\rho)\sigma\\
             &+(\bar{q}_2(t)+\bar{q}_1(t)\rho)\bar{\sigma}\big]\varphi(t)\bigg\}dt+\bar{q}_1(t)dW_1(t)+\bar{q}_2(t)dW_2(t).
\end{split}
\end{equation}

Equating the coefficients of (\ref{bar p1}) and (\ref{bar p2}), gives
\begin{equation}
(\bar{p}(t),\bar{q}_1(t),\bar{q}_2(t))=(-Q(t)\bar{x}(t)-\varphi(t),-\sigma Q(t),-\bar{\sigma}Q(t)),
\end{equation}
where $Q(\cdot)$ is the solution rto the Riccati equation
\begin{equation}\label{Riccati}
\left\{
\begin{aligned}
\begin{split}
&\dot{Q}(t)+\Biggl\{\frac{-\gamma(\sigma^2+\bar{\sigma}^2)^2+2\gamma^2\rho\sigma\bar{\sigma}(\sigma^2+\bar{\sigma}^2)}{(\gamma-1)(\sigma^2+\bar{\sigma}^2)-2\rho\sigma\bar{\sigma}}
 -\frac{4\gamma\rho\sigma\bar{\sigma}(\sigma^2+\bar{\sigma}^2+\rho\sigma\bar{\sigma})}{(\gamma-1)(\sigma^2+\bar{\sigma}^2)-2\rho\sigma\bar{\sigma}}\Biggr\}Q(t)^2\\
&+\frac{2c(\sigma^2+\bar{\sigma}^2+2\rho\sigma\bar{\sigma})}{(\gamma-1)(\sigma^2+\bar{\sigma}^2)-2\rho\sigma\bar{\sigma}}Q(t)-\frac{c^2}{(\gamma-1)(\sigma^2+\bar{\sigma}^2)-2\rho\sigma\bar{\sigma}}=0,\ Q(T)=0,
\end{split}
\end{aligned}
\right.
\end{equation}
and $\varphi(\cdot)$ is a solution to the following equation
\begin{equation}\label{ODE}
\left\{
\begin{aligned}
\begin{split}
&\dot{\varphi}(t)+\bigg\{\bigg[\frac{-\gamma(\sigma^2+\bar{\sigma}^2)^2+2\gamma^2\rho\sigma\bar{\sigma}(\sigma^2+\bar{\sigma}^2)}{(\gamma-1)(\sigma^2+\bar{\sigma}^2)
 -2\rho\sigma\bar{\sigma}}-\frac{4\gamma\rho\sigma\bar{\sigma}(\sigma^2+\bar{\sigma}^2+\rho\sigma\bar{\sigma})}{(\gamma-1)(\sigma^2+\bar{\sigma}^2)-2\rho\sigma\bar{\sigma}}\bigg]Q(t)\\
&+\frac{c(\sigma^2+\bar{\sigma}^2+2\rho\sigma\bar{\sigma})}{(\gamma-1)(\sigma^2+\bar{\sigma}^2)-2\rho\sigma\bar{\sigma}}\bigg\}\varphi(t)
 -\frac{\gamma(\sigma^2+\bar{\sigma}^2)[\frac{1}{2}(\sigma^2+\bar{\sigma}^2+2\rho\sigma\bar{\sigma})+m-r]}{(\gamma-1)(\sigma^2+\bar{\sigma}^2)-2\rho\sigma\bar{\sigma}}Q(t)\\
&+\frac{c[\frac{1}{2}(\sigma^2+\bar{\sigma}^2+2\rho\sigma\bar{\sigma})+m-r]}{(\gamma-1)(\sigma^2+\bar{\sigma}^2)-2\rho\sigma\bar{\sigma}}=0,\ \varphi(T)=0.
\end{split}
\end{aligned}
\right.
\end{equation}

Then, by (\ref{optimal control}) and (\ref{relation}), we can get the optimal control in the following state feedback form:
\begin{equation}\label{feedback}
\begin{split}
\bar{u}(t)=&\ \frac{-\gamma(\sigma^{2}+\bar{\sigma}^{2})Q(t)+c}{(\gamma-1)(\sigma^{2}+\bar{\sigma}^{2})-2\rho\sigma\bar{\sigma}}\bar{x}(t)
-\frac{\gamma(\sigma^{2}+\bar{\sigma}^{2})}{(\gamma-1)(\sigma^{2}+\bar{\sigma}^{2})-2\rho\sigma\bar{\sigma}}\varphi(t)\\&
-\frac{\frac{1}{2}(\sigma^{2}+\bar{\sigma}^{2}+2\rho\sigma\bar{\sigma})+m-r}{(\gamma-1)(\sigma^{2}+\bar{\sigma}^{2})-2\rho\sigma\bar{\sigma}}.
\end{split}
\end{equation}

Next, we apply the approach in Shi \cite{ShiML2005} (see also Shi and Wu \cite{ShiWu2012}), to give the analytical solution to the Riccati equation (\ref{Riccati}). For this target, we rewrite (\ref{Riccati}) as follows
\begin{equation}\label{Riccati2}
\left\{
\begin{aligned}
&\dot{Q}(t)-K_0Q^2(t)+2K_1Q(t)+H=0,\\
&Q(T)=0,
\end{aligned}
\right.
\end{equation}
where we denote
\begin{equation}\label{K0K1H}
\left\{
\begin{split}
K_0=&\frac{\gamma(\sigma^2+\bar{\sigma}^2)^2-2\gamma^2\rho\sigma\bar{\sigma}(\sigma^2+\bar{\sigma}^2)}{(\gamma-1)(\sigma^2+\bar{\sigma}^2)-2\rho\sigma\bar{\sigma}}+\frac{4\gamma\rho\sigma\bar{\sigma}(\sigma^2
     +\bar{\sigma}^2+\rho\sigma\bar{\sigma})}{(\gamma-1)(\sigma^2+\bar{\sigma}^2)-2\rho\sigma\bar{\sigma}},\\
K_1=&\frac{c(\sigma^2+\bar{\sigma}^2+2\rho\sigma\bar{\sigma})}{(\gamma-1)(\sigma^2+\bar{\sigma}^2)-2\rho\sigma\bar{\sigma}},\
     H=-\frac{c^2}{(\gamma-1)(\sigma^2+\bar{\sigma}^2)-2\rho\sigma\bar{\sigma}}.
\end{split}
\right.
\end{equation}

Assuming $\Delta\equiv4(K_1^2+HK_0)>0$, that is,
\begin{equation}\label{suppose}
\frac{\gamma-1}{\gamma^2}<\frac{2\rho\sigma\bar{\sigma}(\sigma^2+\bar{\sigma}^2)}{(\sigma^2+\bar{\sigma}^2)^2+4\rho\sigma\bar{\sigma}(\sigma^2+\bar{\sigma}^2+\rho\sigma\bar{\sigma})},
\end{equation}
we can obtain
\begin{equation}\label{Q}
Q(t)=\frac{\alpha_1+L\alpha_2e^{-\sqrt{\Delta}(T-t)}}{1+Le^{-\sqrt{\Delta}(T-t)}},
\end{equation}
where
\begin{equation}
L\equiv\frac{K_1+\frac{\sqrt{\Delta}}{2}}{K_1-\frac{\sqrt{\Delta}}{2}},\ \alpha_1\equiv\frac{K_1+\frac{\sqrt{\Delta}}{2}}{K_0},\ \alpha_2\equiv\frac{K_1-\frac{\sqrt{\Delta}}{2}}{K_0}.
\end{equation}

Using the denotation (\ref{K0K1H}) we rewrite equation (\ref{ODE}) as
\begin{equation}\label{ODE1}
\begin{aligned}
\dot{\varphi}(t)+f_1(t)\varphi(t)+f_2(t)=0,\ \varphi(T)=0,
\end{aligned}
\end{equation}
where
\begin{equation}\label{f1f2}
\left\{
\begin{split}
f_1(t)=&-K_0Q(t)+K_1,\\
f_2(t)=&\frac{\frac{1}{2}(\sigma^2+\bar{\sigma}^2+2\rho\sigma\bar{\sigma})+m-r}{(\gamma-1)(\sigma^2+\bar{\sigma}^2)-2\rho\sigma\bar{\sigma}}\big[-\gamma(\sigma^2+\bar{\sigma}^2)Q(t)+c\big].
\end{split}
\right.
\end{equation}
The explicit solution to (\ref{ODE1}) is
\begin{equation}\label{f}
\varphi(t)=\exp\Big\{\int_t^Tf_1(s)ds\Big\}\int_t^Tf_2(s)\exp\Big\{-\int_s^Tf_1(r)dr\Big\}ds.
\end{equation}

To summarize, we have the following result.

\vspace{1mm}

{\bf Theorem 4.1}\quad {\it Let (\ref{suppose}) holds, then equations (\ref{Riccati}) and (\ref{ODE}) admit unique solutions $Q(\cdot)$ and $\varphi(\cdot)$, and the optimal control of our risk-sensitive stochastic control problem (\ref{state equation})-(\ref{cost functional}) has the feedback form (\ref{feedback}).}

\vspace{1mm}

{\bf Remark 4.1} Considering the case of $0<\rho<1$, by (\ref{suppose}), we obtain
\begin{equation}\label{gamma}
0<\gamma<\frac{2\rho\sigma\bar{\sigma}}{\sigma^2+\bar{\sigma}^2}+1.
\end{equation}
Obviously, compared with the risk-seeking investors in the most portfolio problem without relevant noise, the existence of correlated noise makes the upper bound larger. It is equivalent to expanding the risk to a certain extent, making some investors with low risk seeking become risk-averse.

In the case of $-1<\rho<0$, (\ref{gamma}) can also be obtained from (\ref{suppose}). Obviously, compared with the risk-seeking investors in the most portfolio problem without relevant noise, the existence of correlated noise makes the upper bound smaller. It is equivalent to reducing the risk to a certain extent, making some investors with risk-aversion become low risk seeking.
\vspace{1mm}

{\bf Remark 4.2} Due to (\ref{gamma}), in this paper we could only obtain the explicit solution to our optimal investment problem in the risk-seeking case, since it relies on the solvability of the Riccati equation (\ref{Riccati}). However, as mentioned in the above remark, because of the introduction of the correlated coefficient $\rho$, we could encounter some cases when the investors are risk-averse. Note this was impossible in Shi and Wu \cite{ShiWu2012}.

\section{Numerical Example and Simulation Result}

In this section, we give a simulating numerical example to show how the corresponding optimal investment statregy changes with respect to the risk-sensitive parameter $\gamma$ and correlation coefficient $\rho$, respectively. We first take some parameters depending on the situation of the real market. Then figures are drawn to illustrate some reasonable analysis adapted to the practical situation. In the example below, we make $T=1$, $c=1$, $m=0.55$, $r=0.05$, $\sigma=0.5$ and $\bar{\sigma}=0.3$.

\subsection{Figures}

Before we make an explicit analysis, we have to give the economical explanation of the correlation coefficient $\rho$ and risk-sensitive parameter $\gamma$. For the correlation coefficient $\rho$, we can understand the practical meaning by considering a case in the financial market. One stock may be influenced by two random factors, when $\rho$ increases, which means that the two factors have strong correlation. When one of which is favorable to the rise of the price of the stock, then the higher the favorable probability of another factor becomes, which is equal to increase the gain and loss in certain degree, and improve the risk. On the contrary, if the two factors have weak correlation, even though one factor has favorable impact, another one has opposite effect to the stock, which is equal to neutralize the volatility of the price, and leads to the drop of the risk. For the risk-sensitive parameter $\gamma$, obviously, when $\gamma\in(0,1)$, it means that the investor is risk-seeking. Therefore, when $\gamma$ approaches to zero, the risk-seeking degree is comparatively high, but when $\gamma$ approaches to one, the risk-seeking attitude of the investor is not very strong.

$(i)$ When $t=0$, then $\bar{x}(0)=0.55$.

\begin{figure}[H]
  \centering
  \includegraphics[width=4.5in,height=2.7in]{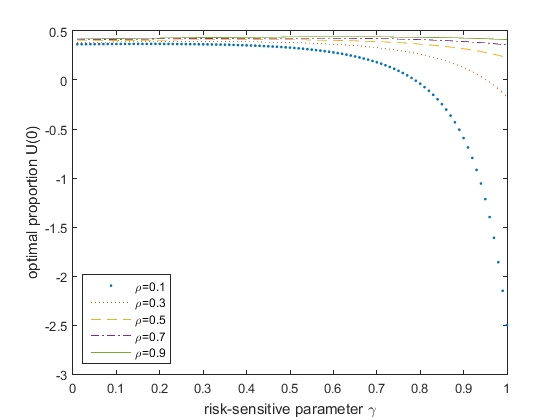}
  \caption{The influence of the risk-sensitive parameter $\gamma$ on the optimal proportion at $t=0$.}
  \label{fig1}
\end{figure}

In the figure 1, firstly, horizontally, it is obviously that the investment proportion of investor has a decreasing trend along with the increase of $\gamma$, where the risk-seeking degree is reducing, which is satisfied with the real situation. Especially, vertically, when $\gamma\in(0,0.4)$, we note that the correlation coefficients $\rho$ have less impact on the investor's choice. More precisely, the investors who have comparatively high risk-seeking attitude can not be easily to change their optimal investment strategies with the change of $\rho$. However, when $\gamma$ is over 0.4, and has a rise, the investors who have comparatively low risk-seeking attitude will become sensitive to the change of $\rho$, that is, correlation coefficient $\rho$ will influence the investors' decision.

\begin{figure}[H]
  \centering
  \includegraphics[width=4.5in,height=2.7in]{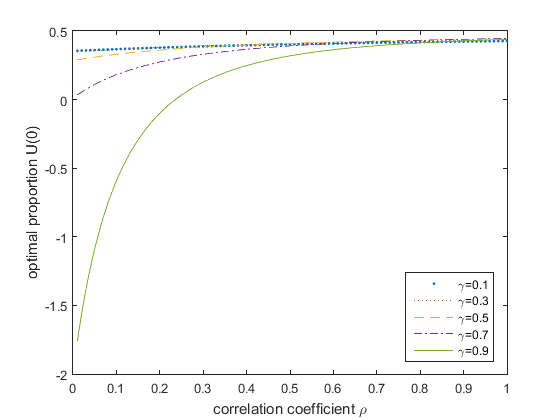}
  \caption{The influence of the correlation coefficient $\rho$ on the optimal proportion at $t=0$.}
  \label{fig2}
\end{figure}

In the figure 2, there is an ascending tendency in the optimal investment proportion following the increase of the correlation coefficient $\rho$.  But an interesting phenomenon is noted, when focusing on the line of $\gamma=0.9$, at first, the optimal investment proportion is a negative value, that is, the investor is short-selling the stock which is different from the others who keep holding the stock. We think the main reason is due to the speciality of the value of $\gamma$ which is very close to 1.

$(ii)$ When $t=0.5$, we define $k_{0.5}=ln\frac{S_{1}(0.5)}{S_{1}(0)}$, which is the log of the ratio of the stock price at time $t=0.5$ to $t=0$. So $\bar{x}(0.5)=k_{0.5}+0.275$.

\begin{figure}[H]
  \centering
  \includegraphics[width=4.5in,height=2.7in]{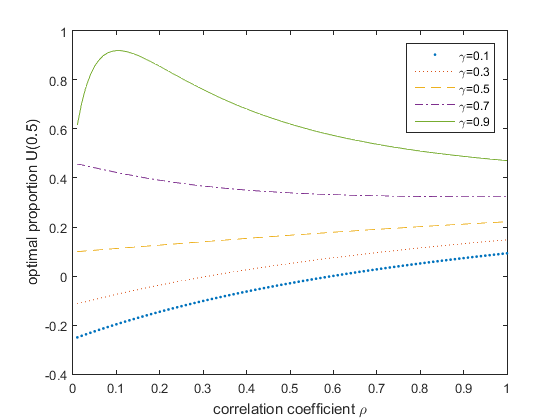}
  \caption{The influence of the correlation coefficient $\rho$ on the optimal proportion at $t=0.5$.}
  \label{fig3}
\end{figure}

When taking $k_{0.5}=0.5$, we get the figure 3. In this case, there exists an interesting phenomenon that, vertically, the investors who have lower risk-seeking attitude prefers to buy more stocks because of the small fluctuation of the price. But the investors who have higher risk-seeking attitude can not stand the comparatively small fluctuation, so they hold the lowest stock proportion. But horizontally, due to the rise of $\rho$, and the increase of the risk, then investors who have low risk-seeking degree cut down the proportion of investing the stock. Especially, for the line of $\gamma=0.9$, when $\rho\in(0,0.1)$, that is, the risk degree is considered to be small, so the investors choose to buy the stock. But when $\rho$ is over 0.1, and even larger, the investors who has low risk-seeking attitude will decrease the investment proportion of stock. While the ones who have high risk-seeking degree begin to hold the increasingly proportion in the stock, which agrees with the real market.

\begin{figure}[H]
  \centering
  \includegraphics[width=4.5in,height=2.7in]{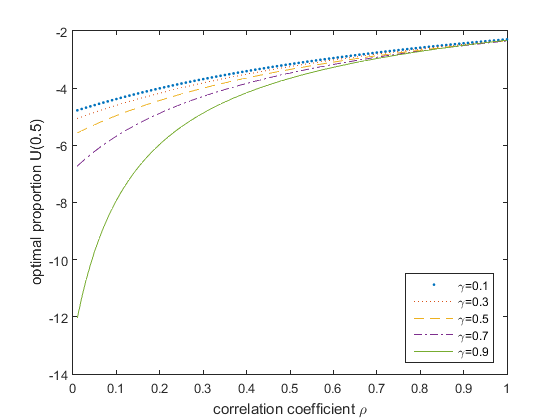}
  \caption{The influence of the correlation coefficient $\rho$ on the optimal proportion at $t=0.5$.}
  \label{fig4}
\end{figure}

When taking $k_{0.5}=2$, we get the figure 4. In this case, apparently, there are a fierce fluctuation in the price of the stock.  Due to the acute rise of the price, it is reasonable for investors to choose to sell the stock. Vertically, the investors who have high risk-seeking degree sell the smallest proportion, which suggests to prefer to keep holding their own stock. Horizontally, the proportion of selling the stock become decreasing along with the increase of $\rho$, that is, although increasing the risk, we are situated in the high profitable situation knowing the increase of the price, the large correlation coefficient $\rho$ enlarge the profitable condition, then investor prefer to keep holding in a long time. All these analysis are fit to the practical situation.

\section{Conclusion}

In this paper, we have introduced a new optimal investment model with correlated random noises in the financial market, which can be reformulated as a kind of a risk-sensitive stochastic control problem with correlated Brownian noises. By using the similar technique in Peng \cite{Peng1990} and Lim and Zhou \cite{LimZhou2005}, necessary and sufficient conditions for the optimal control are proved, where the new first- and second-order adjoint equations and the maximum condition depend on both risk-sensitive parameter and the correlation coefficient. The theoretic results are used to solve the the optimal investment problem with correlated random noises, and the optimal investment strategy is obtained in its state feedback form, under the risk-seeking case. Finally, the sensitivity of the optimal investment strategy to the risk-sensitive parameter and the correlation coefficient, is shown by the numerical simulation and some figures.

Possible extension to the problems with constraints and transaction costs, partial information, jump-diffusion, etc., are interesting topics. We will consider them in our future works.

\appendix

\section{Appendix}

This appendix is devoted to prove Theorem 3.1 and Theorem 3.2. First, we consider the following stochastic control problem:
\begin{align}\label{problem2}
\begin{cases}
{\rm Maximize}:\quad J^\theta(s,x,y;u(\cdot))=\mathbb{E}\big[\exp\big\{\theta[g(x(T))+y(T)]\big\}\big],\\
{\rm subject\ to}:\quad u(\cdot)\in\mathcal{U}[s,T],\\
dx(t)=b(t,x(t),u(t))dt+\sigma_1(t,x(t),u(t))dW_1(t)\\
\hspace{1.4cm}+\sigma_2(t,x(t),u(t))dW_2(t),\ x(s)=x, \\
dy(t)=f(t,x(t),u(t))dt,\ y(s)=y.
\end{cases}
\end{align}
Obviously, (\ref{problem}) reduces to the case of $y=0$ in (\ref{problem2}). The value function $v^\theta:[0,T]\times\mathbb{R}^n\times\mathbb{R}\rightarrow\mathbb{R}$ is
\begin{equation}\label{value function2}
v^\theta(s,x,y):=\sup_{u\in\mathcal{U}[s,T]}J^\theta(s,x,y;u(\cdot)).
\end{equation}
Note that $v^\theta(s,x,y):=\exp\{\theta y\}v^\theta(s,x)$, where $v^\theta(s,x)$ is defined by (\ref{value function}). And {\bf(B3)} implies that $v^\theta\in C^{1,3,\infty}([0,T]\times\mathbb{R}^n\times\mathbb{R})$.

Suppose $(\bar{x}(\cdot),\bar{y}(\cdot),\bar{u}(\cdot))$ is an optimal triple. Assuming {\bf(B1), (B2)}, we can apply the similar technique of Peng \cite{Peng1990} to (\ref{problem2}) to obtain the following first- and second-order adjoint equations:
\begin{equation}\label{adjoint1}
\left\{
\begin{aligned}
&dp(t)=-\Bigg\{\begin{bmatrix}
               \bar{b}_x(t)&0\\
               \bar{f}_x(t)&0
               \end{bmatrix}^\top p(t)+\begin{bmatrix}
               \bar{\sigma}_{1x}(t)&0\\
               0&0
               \end{bmatrix}q_1(t)+\begin{bmatrix}
               \bar{\sigma}_{2x}(t)&0\\
               0&0
               \end{bmatrix}q_2(t)\Bigg\}dt\\
&\qquad\qquad+q_1(t)dW_1(t)+q_2(t)dW_2(t),\\
&p(T)=-\theta\exp\big\{\theta[g(\bar{x}(T)+\bar{y}(T))]\big\}\begin{bmatrix}
               g_x(\bar{x}(T))\\
               1
               \end{bmatrix},
\end{aligned}
\right.
\end{equation}
\begin{equation}\label{adjoint2}
\left\{
\begin{aligned}
&dP(t)=-\Bigg\{\begin{bmatrix}
               \bar{b}_x(t)&0\\
               \bar{f}_x(t)&0
               \end{bmatrix}^\top P(t)+P(t)\begin{bmatrix}
               \bar{b}_x(t)&0\\
               \bar{f}_x(t)&0
               \end{bmatrix}+\begin{bmatrix}
               \bar{\sigma}_{1x}(t)&0\\
               0&0
               \end{bmatrix}^\top P(t)\begin{bmatrix}
               \bar{\sigma}_{1x}(t)&0\\
               0&0
               \end{bmatrix}\\
&\qquad\qquad\quad+\begin{bmatrix}
              \bar{\sigma}_{2x}(t)&0\\
               0&0
              \end{bmatrix}^\top P(t)\begin{bmatrix}
               \bar{\sigma}_{2x}(t)&0\\
               0&0
              \end{bmatrix}+\begin{bmatrix}
              \bar{\sigma}_{1x}(t)&0\\
              0&0
              \end{bmatrix}^\top Q_1(t)\\
&\qquad\qquad\quad+Q_1(t)\begin{bmatrix}
              \bar{\sigma}_{1x}(t)&0\\
              0&0
              \end{bmatrix}+\begin{bmatrix}
              \bar{\sigma}_{2x}(t)&0\\
              0&0
              \end{bmatrix}^\top Q_2(t)+Q_2(t)\begin{bmatrix}
              \bar{\sigma}_{2x}(t)&0\\
              0&0
              \end{bmatrix}\\
&\qquad\qquad\quad+\begin{bmatrix}
              H_{xx}^\theta(t,\bar{x}(t),\bar{u}(t),p(t),q_1(t),q_2(t))&0\\
              0&0
              \end{bmatrix}\Bigg\}dt\\
&\qquad\qquad+Q_1(t)dW_1(t)+Q_2(t)dW_2(t),\\
&P(T)=-\theta\exp\big\{\theta[g(\bar{x}(T)+\bar{y}(T))]\big\}
\begin{bmatrix}
\theta g_x(\bar{x}(T))g_x(\bar{x}(T))^\top+g_{xx}(\bar{x}(T))&\theta g_x(\bar{x}(T))\\
\theta g_x(\bar{x}(T))^\top&\theta
\end{bmatrix},
\end{aligned}
\right.
\end{equation}
where the Hamiltonian $H^\theta:\mathbb{R}\times \mathbb{R}^{n}\times U\times \mathbb{R}^{n+1}\times \mathbb{R}^{n+1}\rightarrow \mathbb{R}$ is given by
\begin{equation}\label{H3}
\begin{split}
&H^\theta(t,x,u)=\left\langle p,\begin{pmatrix}
             b(t,x,u)\\
             f(t,x,u)
            \end{pmatrix}\right\rangle+\left\langle q_1,\begin{pmatrix}
            \sigma_1(t,x,u)\\
             0
            \end{pmatrix}\right\rangle+\left\langle q_2,
            \begin{pmatrix}
            \sigma_2(t,x,u)\\
             0
            \end{pmatrix}\right\rangle.
\end{split}
\end{equation}
The adjoint equations (\ref{adjoint1}) and (\ref{adjoint2}) are linear BSDEs. Under {\bf(B1), (B2)}, for every admissible triple $(\bar{x}(\cdot),\bar{y}(\cdot),\bar{u}(\cdot))$, there are unique solutions $(\bar{p}(\cdot),\bar{q}_1(\cdot),\bar{q}_2(\cdot))\in L^2_\mathcal{F}(s,T;\mathbb{R}^n)\times L^2_\mathcal{F}(s,T;\mathbb{R}^n)\times L^2_\mathcal{F}(s,T;\mathbb{R}^n)$ and $(\bar{P}(\cdot),\bar{Q}_1(\cdot),\bar{Q}_2(\cdot))\in L^2_\mathcal{F}(s,T;\mathbb{R}^{n\times n})\times L^2_\mathcal{F}(s,T;\mathbb{R}^{n\times n})\times L^2_\mathcal{F}(s,T;\mathbb{R}^{n\times n})$ to (\ref{adjoint1}) and (\ref{adjoint2}), respectively.

The $\mathcal{H}$-function for problem (\ref{problem2}) associated with $(\bar{x}(\cdot),\bar{y}(\cdot),\bar{u}(\cdot))$ is defined by
\begin{equation}\label{H4}
\begin{split}
&\mathcal{H}^\theta(t,x,u)=H^\theta(t,x,u,p(t),q_1(t),q_2(t))-\frac{1}{2}\begin{bmatrix}
    \sigma_1(t,\bar{x}(t),\bar{u}(t))\\
    0
    \end{bmatrix}^\top P(t)\begin{bmatrix}
    \sigma_1(t,\bar{x}(t),\bar{u}(t))\\
    0
    \end{bmatrix}\\
&-\frac{1}{2}\begin{bmatrix}
    \sigma_2(t,\bar{x}(t),\bar{u}(t))\\
    0
    \end{bmatrix}^\top P(t)\begin{bmatrix}
    \sigma_2(t,\bar{x}(t),\bar{u}(t))\\
    0
    \end{bmatrix}\\
&+\frac{1}{2}\begin{bmatrix}
    \sigma_1(t,x(t),u(t))-\sigma_1(t,\bar{x}(t),\bar{u}(t)))\\
    0
    \end{bmatrix}^\top P(t)\begin{bmatrix}
    \sigma_1(t,x(t),u(t))-\sigma_1(t,\bar{x}(t),\bar{u}(t)))\\
    0
    \end{bmatrix}\\
&+\frac{1}{2}\begin{bmatrix}
    \sigma_2(t,x(t),u(t))-\sigma_2(t,\bar{x}(t),\bar{u}(t)))\\
    0
    \end{bmatrix}^\top P(t)\begin{bmatrix}
    \sigma_2(t,x(t),u(t))-\sigma_2(t,\bar{x}(t),\bar{u}(t)))\\
    0
    \end{bmatrix}.
\end{split}
\end{equation}

The maximum principle for (\ref{problem2}) can be stated as follows.

\vspace{1mm}

{\bf Proposition A.1} {\it Let {\bf(B1), (B2)} hold. Let $(\bar{x}(\cdot),\bar{y}(\cdot),\bar{u}(\cdot))$ be an optimal triple for the problem (\ref{problem2}). Then there are unique solutions $(p(\cdot),q_1(\cdot),q_2(\cdot))\in L^2_\mathcal{F}(s,T;\mathbb{R}^{n+1})\times L^2_\mathcal{F}(s,T;\mathbb{R}^{n+1})\times L^2_\mathcal{F}(s,T;\mathbb{R}^{n+1})$ and $(P(\cdot),Q_1(\cdot),Q_2(\cdot))\in L^2_\mathcal{F}(s,T;\mathbb{R}^{(n+1)\times(n+1)})\times L^2_\mathcal{F}(s,T;\mathbb{R}^{(n+1)\times(n+1)})\\\times L^2_\mathcal{F}(s,T;\mathbb{R}^{(n+1)\times(n+1)})$ to (\ref{adjoint1}) and (\ref{adjoint2}), respectively, such that
\begin{equation}\label{mp3}
\begin{split}
&H^\theta(t,\bar{x}(t),\bar{u}(t),p(t),q_1(t),q_2(t))-H^\theta(t,\bar{x}(t),u,p(t),q_1(t),q_2(t))\\
&-\frac{1}{2}\begin{bmatrix}
  \sigma_1(t,x(t),u(t))-\sigma_1(t,\bar{x}(t),\bar{u}(t))\\
  0
  \end{bmatrix}^\top P(t)\begin{bmatrix}
  \sigma_1(t,x(t),u(t))-\sigma_1(t,\bar{x}(t),\bar{u}(t))\\
  0
  \end{bmatrix}\\
&-\frac{1}{2}\begin{bmatrix}
  \sigma_2(t,x(t),u(t))-\sigma_2(t,\bar{x}(t),\bar{u}(t))\\
  0
  \end{bmatrix}^\top P(t)\begin{bmatrix}
  \sigma_2(t,x(t),u(t))-\sigma_2(t,\bar{x}(t),\bar{u}(t))\\
  0
  \end{bmatrix}\\
&\leq0,\ \forall u\in\mathbb{R}^k, \ a.e.t\in[0,T], \ \mathbb{P}\mbox{-}{\rm a.s.},
\end{split}
\end{equation}
or equivalently,
\begin{equation}\label{mp4}
\begin{split}
\mathcal{H}^\theta(t,\bar{x}(t),\bar{u}(t))=\min_{u\in\mathbb{R}^k}\mathcal{H}^\theta(t,\bar{x}(t),u),\ a.e. t\in[s,T],\ \mathbb{P}\mbox{-}{\rm a.s}.
\end{split}
\end{equation}}

Sufficient conditions for the optimality of $(\bar{x}(\cdot),\bar{y}(\cdot),\bar{u}(\cdot))$ are as follows.

\vspace{1mm}

{\bf Proposition A.2} {\it Let {\bf(B1), (B2), (B4)} hold. Let $(\bar{x}(\cdot),\bar{y}(\cdot),\bar{u}(\cdot))$ be an admissible triple $(p(\cdot),q_1(\cdot),q_2(\cdot))$ and $(P(\cdot),Q_1(\cdot),Q_2(\cdot))$ satisfy (\ref{adjoint1}) and (\ref{adjoint2}). Suppose that $g(\cdot)$ is convex, $H^\theta(t,\cdot,\cdot,p(t),q_1(t),q_2(t),P(t),Q_1(t),Q_2(t))$ is concave for all $t\in[0,T],\ \mathbb{P}\mbox{-}$a.s., and (\ref{mp4}) holds. Then $(\bar{x}(\cdot),\bar{y}(\cdot),\bar{u}(\cdot))$ is an optimal triple for (\ref{problem2}).}

\vspace{1mm}

In the following, we will transform the adjoint variables $(p(\cdot),q_1(\cdot),q_2(\cdot))$ and $(P(\cdot),Q_1(\cdot),Q_2(\cdot))$ in certain ways. Let $(\bar{x}(\cdot),\bar{y}(\cdot),\bar{u}(\cdot))$ be an optimal triple for (\ref{problem2}), and
\begin{equation*}
\begin{split}
(p(\cdot),q_1(\cdot),q_2(\cdot))
&\equiv\begin{pmatrix}
   \begin{bmatrix}
   p_1(\cdot)\\
   p_2(\cdot)
   \end{bmatrix},
   \begin{bmatrix}
   q_{11}(\cdot)\\
   q_{12}(\cdot)
   \end{bmatrix},
   \begin{bmatrix}
   q_{21}(\cdot)\\
   q_{22}(\cdot)
   \end{bmatrix}
   \end{pmatrix}\\
&\in L^2_\mathcal{F}(s,T;\mathbb{R}^{n+1})\times L^2_\mathcal{F}(s,T;\mathbb{R}^{n+1})\times L^2_\mathcal{F}(s,T;\mathbb{R}^{n+1})
\end{split}
\end{equation*}
be associated with first-order adjoint variables satisfying (\ref{adjoint1}), where $p_1(\cdot), q_{11}(\cdot),q_{21}(\cdot)\in L^2_\mathcal{F}(s,T;\mathbb{R}^n)$, $p_2(\cdot), q_{12}(\cdot),q_{22}(\cdot)\in L^2_\mathcal{F}(s,T;\mathbb{R})$, and $v^\theta(s,x,y)$ is the corresponding value function. Under {\bf(B3)}, $v^\theta(s,x,y)\in C^{1,3,\infty}([0,T]\times\mathbb{R}^n\times\mathbb{R})$.

We take the following transformation of the first-order adjoint variable:
\begin{equation}\label{transformation1}
\tilde{p}(t)=\frac{1}{\theta}\frac{p(t)}{v(t)},
\end{equation}
where $v(t):=v^\theta(t,\bar{x}(t),\bar{y}(t))>0$.

Next, we will derive the equation for $\tilde{p}(\cdot)\equiv\begin{bmatrix}\bar{p}(\cdot)\\\tilde{p}_2(\cdot)\end{bmatrix}$ where $\bar{p}(\cdot)$ is $\mathbb{R}^n$-valued and $\tilde{p}_2(\cdot)$ is scalar valued. First, notice that $v(\cdot)$ is the value function of problem (\ref{problem2}) which has no running cost. Hence, it satisfies
\begin{equation}
\left\{
\begin{split}
dv(t)&=-p_1(t)^\top\sigma_1(t,\bar{x}(t),\bar{u}(t))dW_1(t)-p_1(t)^\top\sigma_2(t,\bar{x}(t),\bar{u}(t))dW_2(t),\\
 v(T)&=\exp\big\{\theta[g(\bar{x}(T))+\bar{y}(T)]\big\}.
\end{split}
\right.
\end{equation}
On the other hand, rearranging (\ref{transformation1}), we obtain $p(t)=\theta v(t)\tilde{p}(t)$. Applying It\^{o}'s formula, and assuming that $\tilde{p}(\cdot)$ satisfies an equation of the following form:
\begin{equation}
d\tilde{p}(t)=\alpha(t)dt+\tilde{q}_1(t)dW_1(t)+\tilde{q}_2(t)dW_2(t),
\end{equation}
we obtain
\begin{equation}
\begin{split}
dp(t)&=\theta v(t)d\tilde{p}(t)+\theta\tilde{p}(t)\big[-p_1(t)^\top\sigma_1(t,\bar{x}(t),\bar{u}(t))dW_1(t)-p_1(t)^\top\sigma_2(t,\bar{x}(t),\bar{u}(t))dW_2(t)\big]\\
     &\quad-\big[\theta p_1(t)^\top\sigma_1(t,\bar{x}(t),\bar{u}(t))\tilde{q}_1(t)+\theta p_1(t)^\top\sigma_2(t,\bar{x}(t),\bar{u}(t))\tilde{q}_2(t)\\
     &\quad+\theta\rho p_1(t)^\top\sigma_2(t,\bar{x}(t),\bar{u}(t))\tilde{q}_1(t)+\theta\rho p_1(t)^\top\sigma_1(t,\bar{x}(t),\bar{u}(t))\tilde{q}_2(t)\big]dt.
\end{split}
\end{equation}
Noting that
\begin{equation}
p_1(t)=\theta v(t)\bar{p}(t),
\end{equation}
we obtain the following expression:
\begin{equation}\label{tilde p}
\begin{split}
d\tilde{p}(t)=&\frac{1}{\theta v(t)}dp(t)+\theta\tilde{p}(t)\bar{p}\sigma_1(t,\bar{x}(t),\bar{u}(t))dW_1(t)+\theta\tilde{p}(t)\bar{p}\sigma_2(t,\bar{x}(t),\bar{u}(t))dW_2(t)]\\
              &+[\theta\tilde{q}_1(t)\bar{p}(t)^\top\sigma_1(t,\bar{x}(t),\bar{u}(t))+\theta\tilde{q}_2(t)\bar{p}(t)^\top\sigma_2(t,\bar{x}(t),\bar{u}(t))\\
              &+\theta\rho\tilde{q}_1(t)\bar{p}(t)^\top\sigma_2(t,\bar{x}(t),\bar{u}(t))+\theta\rho\tilde{q}_2(t)\bar{p}(t)^\top\sigma_1(t,\bar{x}(t),\bar{u}(t))]dt.
\end{split}
\end{equation}
Substituting the expression (\ref{adjoint1}) for $dp(t)$ into (\ref{tilde p}), we can find that the diffusion terms of $d\tilde{p}(t)$ are
\begin{equation}\label{q1}
\begin{split}
&\tilde{q}_1(t)\equiv\begin{bmatrix}
  \bar{q}_1(t)\\
  \tilde{q}_{12}(t)
  \end{bmatrix}=\frac{1}{\theta v(t)}q_1(t)+\theta{\bar{p}(t)}^\top\sigma_1(t,\bar{x}(t),\bar{u}(t))\tilde{p}(t),
\end{split}
\end{equation}
\begin{equation}\label{q2}
\begin{split}
&\tilde{q}_2(t)\equiv\begin{bmatrix}
  \bar{q}_2(t)\\
  \tilde{q}_{22}(t)
  \end{bmatrix}=\frac{1}{\theta v(t)}q_2(t)+\theta{\bar{p}(t)}^\top\sigma_2(t,\bar{x}(t),\bar{u}(t))\tilde{p}(t),
\end{split}
\end{equation}
where $\bar{q}_1(\cdot)$, $\bar{q}_2(\cdot)$ are $\mathbb{R}^n$-valued. Substituting (\ref{q1}), (\ref{q2}) back to (\ref{tilde p}) and using (\ref{adjoint1}), it follows that the transformed first-order adjoint variable $\tilde{p}(\cdot)$ satisfied the following equation, where the terminal condition for $\tilde{p}(T)$ is easily determined from (\ref{transformation1}):
\begin{equation}\label{tilde p2}
\left\{
\begin{aligned}
d\tilde{p}(t)&=-\Bigg\{\begin{bmatrix}
  \bar{b}_x(t)&0\\
  \bar{f}_x(t)&0
  \end{bmatrix}^\top\tilde{p}(t)+\begin{bmatrix}
  \bar{\sigma}_{1x}(t)&0\\
  0&0
  \end{bmatrix}^\top[\tilde{q}_1(t)-\theta\bar{p}(t)^\top\bar{\sigma}_1(t)\tilde{p}(t)]\\
&\qquad+\begin{bmatrix}
  \bar{\sigma}_{2x}(t)&0\\
  0&0
  \end{bmatrix}^\top[\tilde{q}_2(t)-\theta\bar{p}(t)^\top\bar{\sigma}_2(t)\tilde{p}(t)]-\theta[\bar{p}(t)^\top\bar{\sigma}_1(t)\tilde{q}_1(t)\\
&\qquad+\rho\bar{p}(t)^\top\bar{\sigma}_2(t)\tilde{q}_1(t)+\rho\bar{p}(t)^\top\bar{\sigma}_1(t)\tilde{q}_2(t)+\bar{p}(t)^\top\bar{\sigma}_2(t)\tilde{q}_2(t)]\Bigg\}dt\\
&\quad+\tilde{q}_1(t)dW_1(t)+\tilde{q}_2(t)dW_2(t),\\
\tilde{p}(T)&=-\begin{bmatrix}
  g_x(\bar{x}(T))\\
  1
  \end{bmatrix}.
\end{aligned}
\right.
\end{equation}
Expanding (\ref{tilde p2}), it can be easily get
\begin{equation}\label{p2}
\tilde{p}_2(t)=-1,\ \tilde{q}_{12}(t)=\tilde{q}_{22}(t)=0,\ \forall t\in[0,T],
\end{equation}
and $(\bar{p}(\cdot),\bar{q}_1(\cdot),\bar{q}_2(\cdot))\in L^2_\mathcal{F}(s,T;\mathbb{R}^n)\times L^2_\mathcal{F}(s,T;\mathbb{R}^n)\times L^2_\mathcal{F}(s,T;\mathbb{R}^n)$ is a solution to (\ref{BSDE1}). This explains how (\ref{BSDE1}) is derived. Since our derivation can be reserved, it follows from the uniqueness property of (\ref{adjoint1}) that this solution is unique.

Finally, since $\tilde{p}_2(t)=-1$ and $p(t)=\theta v(t)\tilde{p}(t)$, the last component of the extended first order adjoint $p_2(t)=-\theta v(t)$ is essentially the value function.

Now, let $(P(\cdot),Q_1(\cdot),Q_2(\cdot))\in L^2_\mathcal{F}(s,T;\mathbb{R}^{(n+1)\times(n+1)})\times L^2_\mathcal{F}(s,T;\mathbb{R}^{(n+1)\times(n+1)})\times L^2_\mathcal{F}(s,T;\\\mathbb{R}^{(n+1)\times(n+1)})$ be the second-order adjoint variables satisfying (\ref{adjoint2}). We propose the following transformation of the second-order adjoint variable:
\begin{equation}\label{transformation4}
\tilde{P}(t)=\frac{P(t)}{\theta v(t)}+\theta\tilde{p}(t)\tilde{p}(t)^\top\equiv\Gamma(t)+\theta\tilde{p}(t)\tilde{p}(t)^\top.
\end{equation}
Assuming that
\begin{equation}\label{XY1Y2}
d\Gamma(t)=X(t)dt+Y_1(t)dW_1(t)+Y_2(t)dW_2(t),
\end{equation}
for some processes $X(\cdot),Y_1(\cdot),Y_2(\cdot)\in L^2_\mathcal{F}(0,T;\mathbb{R}^{(n+1)\times(n+1)})$, it follows from It\^{o}'s formula and (\ref{transformation4}) that
\begin{equation}\label{P}
\begin{split}
d\Bigg\{\frac{P(t)}{\theta v(t)}\Bigg\}&\equiv d\Gamma(t)=\frac{1}{\theta v(t)}dP(t)+\theta\bar{p}(t)^\top\bar{\sigma_1}(t)\Bigg\{\frac{P(t)}{\theta v(t)}\Bigg\}dW_1(t)\\
                                       &\quad+\theta\bar{p}(t)^\top\bar{\sigma_2}(t)\Bigg\{\frac{P(t)}{\theta v(t)}\Bigg\}dW_2(t)+\big[\theta\bar{p}(t)^\top\bar{\sigma_1}(t)Y_1(t)\\
                                       &\quad+\theta\rho\bar{p}(t)^\top\bar{\sigma_1}(t)Y_2(t)+\theta\rho\bar{p}(t)^\top\bar{\sigma_2}(t)Y_1(t)+\theta\bar{p}(t)^\top\bar{\sigma_2}(t)Y_2(t)\big]dt.
\end{split}
\end{equation}
Substituting the expression for $dP(t)$ and notice (\ref{XY1Y2}), gives
\begin{equation}
Y_i(t)=\frac{Q_i(t)}{\theta v(t)}+\theta\bar{p}(t)^\top\bar{\sigma}_i(t)\tilde{P}(t)-\theta^2\bar{p}(t)^\top\bar{\sigma}_i(t)\tilde{p}(t)\tilde{p}(t)^\top,\ i=1,2.
\end{equation}
By (\ref{P}), (\ref{transformation4}) and using It\^{o}'s formula, we obtain
\begin{equation}\label{tilde P}
\begin{split}
d\tilde{P}(t)=&-\Bigg\{\begin{bmatrix}
  \bar{b}_x(t)&0\\
  \bar{f}_x(t)&0
  \end{bmatrix}^\top\tilde{P}(t)+\tilde{P}(t)\begin{bmatrix}
  \bar{b}_x(t)&0\\
  \bar{f}_x(t)&0
  \end{bmatrix}+\begin{bmatrix}
  \bar{\sigma}_{1x}(t)&0\\
  0&0
  \end{bmatrix}^\top\bigg(\frac{P(t)}{\theta v(t)}\bigg)\begin{bmatrix}
  \bar{\sigma}_{1x}(t)&0\\
  0&0
  \end{bmatrix}\\
&\quad+\begin{bmatrix}
  \bar{\sigma}_{2x}(t)&0\\
  0&0
  \end{bmatrix}^\top\bigg(\frac{P(t)}{\theta v(t)}\bigg)\begin{bmatrix}
  \bar{\sigma}_{2x}(t)&0\\
  0&0
  \end{bmatrix}\quad+\begin{bmatrix}
  \bar{\sigma}_{1x}(t)&0\\
  0&0
  \end{bmatrix}^\top\Big[\frac{Q_1(t)}{\theta v(t)}\\
&\quad+\theta\tilde{q}_1(t)\tilde{p}(t)^\top-\theta^2\bar{p}(t)^\top\bar{\sigma}_1(t)\tilde{p}(t)\tilde{p}(t)^\top\Big]+\Big[\frac{Q_1(t)}{\theta v(t)}+\theta\tilde{p}(t)\tilde{q}_1(t)^\top\\
&\quad-\theta^2\bar{p}(t)^\top\bar{\sigma}_1(t)\tilde{p}(t)\tilde{p}(t)^\top\Big]
  \begin{bmatrix}
  \bar{\sigma}_{1x}(t)&0\\
  0&0
  \end{bmatrix}+\begin{bmatrix}
  \bar{\sigma}_{2x}(t)&0\\
  0&0
  \end{bmatrix}^\top\Big[\frac{Q_2(t)}{\theta v(t)}+\theta\tilde{q}_2(t)\tilde{p}(t)^\top\\
&\quad-\theta^2\bar{p}(t)^\top\bar{\sigma}_2(t)\tilde{p}(t)\tilde{p}(t)^\top\Big]
  +\Big[\frac{Q_2(t)}{\theta v(t)}+\theta\tilde{p}(t)\tilde{q}_2(t)^\top-\theta^2\bar{p}(t)^\top\bar{\sigma}_2(t)\tilde{p}(t)\tilde{p}(t)^\top\Big]\\
&\quad\times\begin{bmatrix}
  \bar{\sigma}_{2x}(t)&0\\
  0&0
  \end{bmatrix}+\frac{1}{\theta v(t)}\begin{bmatrix}
  H_{xx}^{\theta}(t,\bar{x}(t),\bar{u}(t),p(t),q_1(t),q_2(t))&0\\
  0&0
  \end{bmatrix}\\
&\quad-\theta\bar{p}(t)^\top\bar{\sigma}_1(t)\{Y_1(t)+\theta\tilde{p}(t)\tilde{q}_1(t)^\top+\theta\tilde{q}_1(t)\tilde{p}(t)^\top\}-\theta\bar{p}(t)^\top\bar{\sigma}_2(t)\{Y_2(t)\\
&\quad+\theta\tilde{p}(t)\tilde{q}_2(t)^\top+\theta\tilde{q}_2(t)\tilde{p}(t)^\top\}-\theta\rho\bar{p}(t)^\top\bar{\sigma}_1(t)Y_2(t)-\theta\rho\bar{p}(t)^\top\bar{\sigma}_2(t)Y_1(t)\\
&\quad-\theta^2\rho\tilde{p}(t)\tilde{q}_1(t)^\top\bar{\sigma}_2(t)^\top\bar{p}(t)-\theta^2\rho\tilde{p}(t)\tilde{q}_2(t)^\top\bar{\sigma}_1(t)^\top\bar{p}(t)-\theta\tilde{q}_1(t)\tilde{q}_1(t)^\top\\
&\quad-\theta^2\rho\bar{p}(t)^\top\bar{\sigma}_2(t)\tilde{q}_1(t)\bar{p}(t)^\top-\theta^2\rho\bar{p}(t)^\top\bar{\sigma}_{1}(t)\tilde{q}_{2}(t)\bar{p}(t)^\top-\theta\rho\tilde{q}_1(t)\tilde{q}_2(t)^\top\\
&\quad-\theta\rho\tilde{q}_2(t)\tilde{q}_1(t)^\top-\theta\tilde{q}_2(t)\tilde{q}_2(t)^\top\Bigg\}dt+\tilde{Q}_1(t)dW_1(t)+\tilde{Q}_2(t)dW_2(t),
\end{split}
\end{equation}
where
\begin{equation}\label{tilde Q1Q2}
\tilde{Q}_i(t)=\frac{Q_i(t)}{\theta v(t)}+\theta\bar{p}(t)^\top\bar{\sigma}_i(t)[\tilde{P}(t)-\theta\tilde{p}(t)\tilde{p}(t)^\top]+\theta\tilde{p}(t)\tilde{q}_i(t)^\top+\theta\tilde{q}_i(t)\tilde{p}(t)^\top,\ i=1,2.
\end{equation}
Reminding the definition of $\bar{H}^\theta$, we have
\begin{equation}\label{vtheta}
\begin{split}
\frac{1}{\theta v(t)}\begin{bmatrix}
  H_{xx}^\theta(t,\bar{x}(t),\bar{u}(t),p(t),q_1(t),q_2(t))&0\\
  0&0
  \end{bmatrix}
=\begin{bmatrix}
  \bar{H}_{xx}^{\theta}(t,\bar{x}(t),\bar{u}(t),p(t),q_1(t),q_2(t))&0\\
  0&0
\end{bmatrix}.
\end{split}
\end{equation}
Noticing that
\begin{equation}
Y_i(t)=\tilde{Q}_i(t)-\theta\tilde{p}(t)\tilde{q}_i(t)^\top-\theta\tilde{q}_i(t)\tilde{p}(t)^\top,\ i=1,2,
\end{equation}
together with the transformation (\ref{q1}), (\ref{q2}), (\ref{transformation4}), (\ref{tilde P}), (\ref{tilde Q1Q2}) and (\ref{vtheta}), we can obtain
\begin{equation*}
\left\{
\begin{aligned}
d\tilde{P}(t)=&-\Bigg\{\begin{bmatrix}
  \bar{b}_x(t)&0\\
  \bar{f}_x(t)&0
  \end{bmatrix}^\top\tilde{P}(t)+\begin{bmatrix}
  \bar{\sigma}_{1x}(t)&0\\
  0&0
  \end{bmatrix}^\top[\tilde{P}(t)-\theta\tilde{p}(t)\tilde{p}(t)^\top]\begin{bmatrix}
  \bar{\sigma}_{1x}(t)&0\\
  0&0
  \end{bmatrix}\\
&\quad+\tilde{P}(t)\begin{bmatrix}
  \bar{b}_{x}(t)&0\\
  \bar{f}_{x}(t)&0
  \end{bmatrix}+\begin{bmatrix}
  \bar{\sigma}_{2x}(t)&0\\
  0&0
  \end{bmatrix}^\top[\tilde{P}(t)-\theta\tilde{p}(t){\tilde{p}(t)}^\top]\begin{bmatrix}
  \bar{\sigma}_{2x}(t)&0\\
  0&0
  \end{bmatrix}\\
&\quad+\begin{bmatrix}
  \bar{\sigma}_{1x}(t)&0\\
  0&0
  \end{bmatrix}^\top[\tilde{Q}_1(t)-\theta\bar{p}(t)^\top\bar{\sigma}_1(t)\tilde{P}(t)-\theta\tilde{p}(t)\tilde{q}_1(t)^\top]\\
&\quad+[\tilde{Q}_1(t)-\theta\bar{p}(t)^\top\bar{\sigma}_1(t)\tilde{P}(t)-\theta\tilde{q}_1(t)\tilde{p}(t)^\top]^\top
  \begin{bmatrix}
  \bar{\sigma}_{1x}(t)&0\\
  0&0
  \end{bmatrix}^\top\\
&\quad+\begin{bmatrix}\bar{\sigma}_{2x}(t)&0\\
  0&0
  \end{bmatrix}^\top\big[\tilde{Q}_2(t)-\theta\bar{p}(t)^\top\bar{\sigma}_2(t)\tilde{P}(t)-\theta\tilde{p}(t)\tilde{q}_2(t)^\top\big]\\
&\quad+\big[\tilde{Q}_2(t)-\theta\bar{p}(t)^\top\bar{\sigma}_2(t)\tilde{P}(t)-\theta\tilde{q}_2(t)\tilde{p}(t)^\top\big]^\top
  \begin{bmatrix}
  \bar{\sigma}_{2x}(t)&0\\
  0&0
  \end{bmatrix}^\top\\
&\quad+\begin{bmatrix}
  \bar{H}_{xx}^{\theta}(t,\bar{x}(t),\bar{u}(t),p(t),q_1(t),q_2(t))&0\\
  0&0
  \end{bmatrix}-\theta\bar{p}(t)^\top\bar{\sigma}_1(t)\tilde{Q}_1(t)-\theta\bar{p}(t)^\top\bar{\sigma}_2(t)\tilde{Q}_2(t)\\
&\quad-\theta\rho\bar{p}(t)^\top\bar{\sigma}_2(t)\tilde{Q}_1(t)-\theta\rho\bar{p}(t)^\top\bar{\sigma}_1(t)\tilde{Q}_2(t)+\theta^2\rho\bar{p}(t)^\top\bar{\sigma}_1(t)\tilde{p}(t)\tilde{q}_2(t)^\top\\
&\quad+\theta^2\rho\bar{p}(t)^\top\bar{\sigma}_1(t)\tilde{q}_2(t)\tilde{p}(t)^\top+\theta^2\rho\bar{p}(t)^\top\bar{\sigma}_2(t)\tilde{p}(t)\tilde{q}_1(t)^\top
 +\theta^2\rho\bar{p}(t)^\top\bar{\sigma}_2(t)\tilde{q}_2(t)\tilde{p}(t)^\top\\
&\quad-\theta^2\rho\tilde{p}(t)\tilde{q}_1(t)^\top\bar{\sigma}_2(t)^\top\bar{p}(t)-\theta^2\rho\tilde{p}(t)\tilde{q}_2(t)^\top\bar{\sigma}_1(t)^\top\bar{p}(t)
-\theta^2\rho\bar{p}(t)^\top\bar{\sigma}_2(t)\tilde{q}_1(t)\tilde{p}(t)^\top\\
&\quad-\theta^2\rho\bar{p}(t)^\top\bar{\sigma}_1(t)\tilde{q}_2(t)\tilde{p}(t)^\top-\theta\tilde{q}_1(t)\tilde{q}_1(t)^\top-\theta\rho\tilde{q}_1(t)\tilde{q}_2(t)^\top-\theta\rho\tilde{q}_2(t)\tilde{q}_1(t)^\top\\
&\quad-\theta\tilde{q}_2(t)\tilde{q}_2(t)^\top\Bigg\}dt+\tilde{Q}_1(t)dW_1(t)+\tilde{Q}_2(t)dW_2(t),\\
\tilde{P}(T)=&-\begin{bmatrix}
  g_{xx}(\bar{x}(T))&0\\
  0&0
  \end{bmatrix}.
\end{aligned}
\right.
\end{equation*}
Therefore, it follows that
\begin{align*}
&\tilde{P}(t)=\begin{bmatrix}
  \bar{P}(t)&0\\
  0&0
  \end{bmatrix},\quad\tilde{Q}_1(t)=\begin{bmatrix}
  \bar{Q}_1(t)&0\\
  0&0
  \end{bmatrix},\quad\tilde{Q}_2(t)=
  \begin{bmatrix}
  \bar{Q}_2(t)&0\\
  0&0
  \end{bmatrix},
\end{align*}
where $(\bar{P}(\cdot),\bar{Q}_1(\cdot),\bar{Q}_2(\cdot))\in L^2_\mathcal{F}(s,T;\mathbb{R}^{n\times n})\times L^2_\mathcal{F}(s,T;\mathbb{R}^{n\times n})\times L^2_\mathcal{F}(s,T;\mathbb{R}^{n\times n})$ is the solution to (\ref{BSDE2}). As in the first-order case, the solution is unique.

Now, we consider (\ref{mp1}). By (\ref{transformation1}), (\ref{q1}), (\ref{q2}), (\ref{p2}), we have
\begin{equation}
\begin{split}
H^\theta(t,x,u,p(t),q_1(t),q_2(t))=[\theta v(t)]\bar{H}^\theta(t,x,u,p(t),q_1(t),q_2(t)),
\end{split}
\end{equation}
where $\bar{H}^{\theta}(t,x,u,,q_{1},q_{2})$ is given by (\ref{H1}), and
\begin{equation*}
\begin{split}
 &\frac{1}{2}\begin{bmatrix}
  \sigma_1(t,\bar{x}(t),\bar{u}(t))-\sigma_1(t,\bar{x}(t),u(t))\\
  0
  \end{bmatrix}^\top P(t)\begin{bmatrix}
  \sigma_1(t,\bar{x}(t),\bar{u}(t))-\sigma_1(t,\bar{x}(t),u(t))\\
  0
  \end{bmatrix}\\
 &+\frac{1}{2}\begin{bmatrix}
  \sigma_2(t,\bar{x}(t),\bar{u}(t))-\sigma_2(t,\bar{x}(t),u(t))\\
  0
  \end{bmatrix}^\top P(t)\begin{bmatrix}
  \sigma_2(t,\bar{x}(t),\bar{u}(t))-\sigma_2(t,\bar{x}(t),u(t))\\
  0
  \end{bmatrix}\\
=&\Big[\frac{\theta v(t)}{2}\Big]\Big\{\big[\sigma_1(t,\bar{x}(t),\bar{u}(t))-\sigma_1(t,\bar{x}(t),u(t))]^\top[\bar{P}(t)-\theta\bar{p}(t)\bar{p}(t)^\top][\sigma_1(t,\bar{x}(t),\bar{u}(t))\\
 &-\sigma_1(t,\bar{x}(t),u(t))]+[\sigma_2(t,\bar{x}(t),\bar{u}(t))-\sigma_2(t,\bar{x}(t),u(t))]^\top[\bar{P}(t)-\theta\bar{p}(t)\bar{p}(t)^\top]\\
 &\times[\sigma_2(t,\bar{x}(t),\bar{u}(t))-\sigma_2(t,\bar{x}(t),u(t))]\Big\}.
\end{split}
\end{equation*}
Since $v(t)>0$,  it follows that the maximum condition (\ref{mp3}) is equivalent to
\begin{equation*}
\begin{split}
&\bar{H}^\theta(t,\bar{x}(t),\bar{u}(t),\bar{p}(t),\bar{q}_1(t),\bar{q}_2(t))-\bar{H}^\theta(t,\bar{x}(t),u,\bar{p}(t),\bar{q}_1(t),\bar{q}_2(t))\\
&-\frac{1}{2}[\sigma_1(t,\bar{x}(t),\bar{u}(t))-\sigma_1(t,\bar{x}(t),u(t))]^\top(\bar{P}(t)-\theta\bar{p}(t)\bar{p}(t)^\top)[\sigma_1(t,\bar{x}(t),\bar{u}(t))-\sigma_1(t,\bar{x}(t),u(t))]\\
&-\frac{1}{2}[\sigma_2(t,\bar{x}(t),\bar{u}(t))-\sigma_2(t,\bar{x}(t),u(t))]^\top(\bar{P}(t)-\theta\bar{p}(t)\bar{p}(t)^\top)[\sigma_2(t,\bar{x}(t),\bar{u}(t))-\sigma_2(t,\bar{x}(t),u(t))]\\
&\leq0,
\end{split}
\end{equation*}
which gives us (\ref{mp1}). The equivalent condition (\ref{mp2}) can be get by direct manipulation. This complete the proofs of Theorems 3.1 and 3.2.

\end{document}